\documentclass[reqno,11pt,draft]{amsart}
\usepackage{amsmath,amssymb,amsfonts,amsxtra}

\usepackage{pgf,tikz}
\usetikzlibrary{arrows}
\usepackage{verbatim}
\usepackage{color}
\usepackage{hyperref}

\theoremstyle{plain}
\newtheorem{thm}{Theorem}[section]
\newtheorem{lem}[thm]{Lemma}
\newtheorem{cor}[thm]{Corollary}
\newtheorem{propo}[thm]{Proposition}

\theoremstyle{definition}
\newtheorem{defn}[thm]{Definition}

\newtheorem{parraf}[thm]{}
\newtheorem{subparraf}{}[thm]
\newtheorem*{ack}{Acknowledgments}

\theoremstyle{remark}
\newtheorem*{rem}{Remark}

\numberwithin{equation}{thm}

\newcommand{\hd}{\widehat{d}}

\newcommand{\hh}{\widehat{h}}

\newcommand{\HA}{\widehat{A}}

\newcommand{\WL}{\widehat{L}}

\newcommand{\WCL}{\widehat{\mathcal{L}}}
\newcommand{\tc}{\widehat{\otimes}}
\newcommand{\om}{\widehat{\Omega}}

\newcommand{\NN}{\mathbb N}

\newcommand{\FS}{\mathfrak S}

\newcommand{\FU}{\mathfrak U}
\newcommand{\FV}{\mathfrak V}
\newcommand{\FW}{\mathfrak W}
\newcommand{\FX}{\mathfrak X}
\newcommand{\FY}{\mathfrak Y}
\newcommand{\FZ}{\mathfrak Z}

\newcommand{\cc}{\mathsf{c}}
\newcommand{\ic}{\vec{\mathsf{c}}}
\newcommand{\sch}{\mathsf {Sch}}
\newcommand{\sfn}{\mathsf {FS}}
\newcommand{\scha}{\mathsf {Sch}_{\mathsf {af}}}
\newcommand{\sfna}{\sfn_{\mathsf {af}}}
\newcommand{\A}{\mathsf{A}}

\newcommand{\C}{\boldsymbol{\mathsf{C}}}
\newcommand{\D}{\boldsymbol{\mathsf{D}}}
\newcommand{\R}{\boldsymbol{\mathsf{R}}}
\newcommand{\LD}{\boldsymbol{\mathsf{L}}}
\newcommand{\RLA}{\boldsymbol{\mathsf{\Lambda}}}
\newcommand{\RGA}{\boldsymbol{\mathsf{\Gamma}}}

\newcommand{\CE}{\mathcal E}
\newcommand{\CF}{\mathcal F}
\newcommand{\CG}{\mathcal G}
\newcommand{\CH}{\mathcal H}
\newcommand{\CI}{\mathcal I}
\newcommand{\CJ}{\mathcal J}
\newcommand{\CK}{\mathcal K}
\newcommand{\CL}{\mathcal L}
\newcommand{\CM}{\mathcal M}
\newcommand{\CN}{\mathcal N}
\newcommand{\CO}{\mathcal O}

\newcommand{\dirlim}[1]{\begin{array}[t]{c} \mathrm{lim}\\[-7.5 pt]
 {\longrightarrow} \\[-7.5 pt] {\scriptstyle {#1}} \end{array}}
 
\newcommand{\invlim}[1]{\begin{array}[t]{c} \mathrm{lim}\\[-7.5 pt]
 {\longleftarrow} \\[-7.5 pt] {\scriptstyle {#1}} \end{array}}

\newcommand{\lto}{\longrightarrow}
\newcommand{\xto}{\xrightarrow}
\newcommand{\inc}{\hookrightarrow}
\newcommand{\tr}{\triangle}

\newcommand{\cf}{{\it cf.}}

\DeclareMathOperator{\h}{H}
\DeclareMathOperator{\spec}{Spec}
\DeclareMathOperator{\spf}{Spf}
\DeclareMathOperator{\Ker}{Ker}
\DeclareMathOperator{\Hom}{Hom}

\DeclareMathOperator{\Dercont}{Dercont}
\DeclareMathOperator{\ga}{\Gamma}
\DeclareMathOperator{\LA}{\Lambda}
\DeclareMathOperator{\shom}{\CH\mathit{om}}
\DeclareMathOperator{\ext}{Ext}

\DeclareMathOperator{\modu}{\text{-}\mathsf{mod}}
\DeclareMathOperator{\Modu}{\mathsf{Mod}}
\DeclareMathOperator{\rshom}{\R\!\shom}
\DeclareMathOperator{\topo}{top}

\begin{document}

\title[Deformation of formal schemes]{Deformation of  formal schemes through local homology}

\author[M. P\'erez Rodr\'{\i}guez]{Marta P\'erez Rodr\'{\i}guez}
\address{Departamento de Matem\'a\-ticas\\
Escola Superior de En\-xe\-\~ne\-r\'{\i}a Inform\'atica\\
Campus de Ourense, Univ. de Vigo\\
E-32004 Ou\-ren\-se, Spain}
\email{martapr@uvigo.es}

\thanks{This work has been partially supported by
European Community under the contract AAG MRTN CT 2003 504917, together
Spain's MEC and E.U.'s FEDER research project  MTM2011-26088 and Xunta de Galicia's
PGIDIT10PXIB207144PR}

\date{July 24, 2015}

\hyphenation{pseu-do}

\begin{abstract} 
We define the cotangent complex of a morphism $f \colon \FX \to \FY$ of locally noetherian formal schemes as an object in the derived category $\D^-(\FX)$ through local homology.  We discuss its basic properties and  establish the basics results of a deformation theory, providing a characterization of smooth and \'etale morphisms. This leads to  simpler lifting results depending on a differential module, for a class of non smooth morphism of usual schemes.  We also give descriptions of the cotangent complex in the case of regular closed immersions and complete intersection morphisms of formal schemes.

\end{abstract}

\maketitle

\tableofcontents

\section*{Introduction}
The aim of this paper is to provide a cotangent complex object and a deformation theory in the category of locally noetherian formal schemes. 
We continue the development of a infinitesimal theory of locally noetherian formal schemes established in \cite{AJP1}, \cite{AJP2} and \cite{P}.

Formal schemes  have become relevant in different areas of mathematics as rigid geometry \cite{ray}, cohomology of singular spaces \cite{ha2}, $p$-adic cohomologies  or stable homotopy theory \cite{st}. In parallel, several authors have continued the foundations on formal schemes of \cite{EGA1} and \cite{EGA31} from a cohomological point of view \cite{y},  \cite{AJL2}, \cite{AJL3}, \cite{LNS}. A common hypothesis in all these works is that they use pseudo finite type condition of morphisms. 
In \cite[\S 10]{EGA1} only finite type (therefore \emph{adic}) morphisms are treated. However,  they do not cover a wide range of relevant morphisms, such as projections of formal discs $\mathbb{Z}[[t]]\times \FX \to \FX$ or completions $X_{/X'} \to X$, that appear frequently in different situations. Pseudo finite type morphisms \cite{AJL2} are (\emph{non adic}) morphisms with a finiteness condition that generalizes finite type morphisms of usual schemes. For instance, finite type morphisms, projections of formal discs and completions are of pseudo finite type.

The study of Grothendieck duality on formal schemes \cite{AJL2}, \cite{AJL3} led to a relevant advance in the study of suitable cohomological coefficients in the derived category $\D(\FX)$ of a locally noetherian formal scheme $\FX$.
However, some cohomological tools still lack, for instance, a suitable cotangent complex in order to get a deformation theory for formal schemes.

At the end of the 60s  Illusie  \cite{I} developed a   cotangent complex and deformation theory\footnote{ We should mention that first versions were given by Grothendieck \cite{G} as a one term complex; Lichtenbaum and Schlessinger  \cite{LS} for commutative rings  as a three term complex;  Berthelot \cite{sga6}  for closed immersions of smooth scheme morphisms and as an object in the derived category with coherent cohomology. Finally, Andr\'e \cite{An} and Quillen \cite{Q} introduced it in a simplicial commutative ring setting.}  in a simplicial topos setting. In the last decade interest have arisen  in   deformation theory on different algebra geometric contexts. Laumon and Moret-Bailly \cite{LMB} define the  cotangent complex of a $1$-morphism of algebraic stacks as an ind-object in the bounded above derived category.  This approach was used  to provide a deformation theory of algebraic stacks by Aoki \cite{Ao}, and in the case of a scheme over an algebraic stack  by  Olsson \cite{O2}. In logarithmic geometry, Olsson introduces the logarithmic cotangent complex as a ind-object in the bounded above derived category \cite{O1}. 
In derived algebraic geometry, To\"en and Vezzosi  \cite{TV} gave a derived interpretation of the cotangent complex, defining it  by representability of the derivations functor and establishing an obstruction theory.  

On the other hand, in almost geometry, Gabber and Ramero  \cite{GR}  analyzed deformation theory through the cotangent complex. They also defined a cotangent complex object for  certain formal schemes and adic spaces. In their approach an important hypothesis is that morphisms of formal schemes are of  locally finite type (\emph{adic}) and over a rank one valued field. They introduce the cotangent complex in an affine context and then define it globally through a sheafification procedure.

In the context of  formal schemes,  the cotangent complex of ringed topoi does not
 have the right properties, since it does not
 consider the topology of the structural sheaves. For instance, with pseudo finite type hypothesis it is not of coherent cohomology. A first approach to obstruction theory of formal schemes was given  in \cite{P} through the complete differential module and under smoothness (therefore pseudo finite type) and separation hypothesis.  This fact and the previous discussion motivates us to provide a general theory of cotangent complex and deformation  for pseudo finite type morphisms of locally noetherian formal schemes.

Even though our treatment generalizes the deformation theory for  schemes (\cite[Chapitres II and III]{I}),  it does not follow from it. Our arguments are based on properties of the \emph{complete} differential module \cite{AJP1}, the derived category of complexes with coherent cohomology associated to a formal scheme and the homology localization functor \cite{AJL2} and Greenlees-May duality  \cite{AJL3}. We use basics of the  theory of the cotangent complex and deformations on ringed topoi  \cite{I},   considering the cotangent complex of ringed spaces in the bounded above derived category  via Dold-Puppe equivalence \cite{DP}.
We expect that our results may be applied to the  cohomological study of singular varieties.

\vspace{10pt}

Let us now describe  the contents of this paper.
In Section \ref{sec1} we summarize the main background material about formal schemes that will be use along the exposition: pseudo-finite type morphisms, infinitesimal conditions and differential pair $(\om^1_{\FX/\FY}, \hd_{\FX/\FY})$, among others.

In Section \ref{sec2} we introduce the cotangent complex   in an affine formal geometry context.   We define the \emph{complete cotangent complex} $\WL_{B/A} \in \C^{\le 0}(B)$ of a continuous map of adic rings $A \to B$ and show that basic properties, as the existence of an augmentation map and functoriality, hold  in $ \C^{\le 0}(B)$. The \emph{pseudo finite type} hypothesis allows us to establish that  $\WL_{B/A}$ is a complex with coherent bounded above cohomology such that $\h^0(\WL_{B/A})= \om^1_{B/A}$. Moreover, given  another morphism $B \to C$ there is a distinguished triangle  $\WL_{B/A} \otimes_{B}C  \to \WL_{C /A} \to \WL_{C /B} \overset{+} \to$ and we obtain specific descriptions of $\WL_{B/A}$ for smooth and \'etale morphisms and closed immersions. The results of this section will be essential in what follows.

Section \ref{sec3} is dedicated to the study of the cotangent complex object of a  morphism $\FX \to \FY$ of locally noetherian formal schemes. We introduce the complete cotangent complex as $\WCL_{\FX/\FY}:= \RLA_{\FX}(\CL_{\FX/\FY})  \in \D^-(\FX)$, where $\RLA_{\FX}$ is the homology localization functor and $\CL_{\FX/\FY}$ is the cotangent complex of ringed spaces. We prove that a pseudo finite type map has  the  properties considered desirable for a cotangent complex object, namely,
\begin{itemize}
\item
\textbf{augmentation map} (\ref{sheafaug}):
$\WCL_{\FX/\FY} \to  \om^1_{\FX/\FY}$;
\item
\textbf{functoriality} (\ref{functoriality}, Proposition \ref{coh}): given  a morphism $\FY \to \FS$
there is a morphism in $\D_\cc^-(\FX)$ 
\[
\LD f^{*} \WCL_{\FY/\FS} \to \WCL_{\FX/\FS}
\]
compatible with the augmentation maps; 
 \item
\textbf{localization} (Proposition \ref{localcotag}): if $\FU \overset{i}\inc \FX$ is an open immersion, then $i^*\WCL_{\FX/\FY} = \WCL_{\FU/\FY}$;
\item
\textbf{coherence} (Proposition \ref{coh}, Corollary \ref{h0}):
$\WCL_{\FX/\FY}  \in \D_\cc^-(\FX)$, being $\h^{0}(\WCL_{\FX/\FY}) \simeq \om^1_{\FX/\FY}$ and $\h^{i}(\WCL_{\FX/\FY}) \simeq(\h^{i}(\WL_{B/A}))^{\tr}$ in the affine case;
\item
\textbf{distinguished triangle} (Proposition \ref{triancoh}): there is a distinguished triangle in $ \D_\cc^-(\FX)$ 
\[
\LD f^{*} \WCL_{\FY/\FS} \to \WCL_{\FX/\FS} \to \WCL_{\FX/\FY} \xto{+}; 
\]
\item
\textbf{flat base-change} (Proposition \ref{flatchange}): given $\FX$ or $\FY'$ flat over $\FY$  and $u \colon \FX':= \FX \times_{\FY} \FY'\to \FX$, then $\LD u^* \WCL_{\FX/\FY} \simeq \WCL_{\FX'/\FY'}$.
\end{itemize}
\vspace{5pt}

Deformation theory for locally noetherian formal schemes is treated in Section \ref{sec4}. The key results are Proposition \ref{extenfs}, which ensures that an extension of a formal scheme  by a square zero coherent module is a formal scheme, and the isomorphisms 
\begin{alignat*}{2}
\ext^{i}_{\CO_{\FX}} (\WCL_{\FX/\FY}, \CF) &\simeq \ext^{i}_{\CO_{\FX}} (\CL_{\FX/\FY}, \CF)\\
\ext^{i}_{\CO_{\FX}} (\LD f^*\WCL_{\FY/\FS}, \CF) &\simeq \ext^{i}_{\CO_{\FX}} (f^*\CL_{\FY/\FS}, \CF)
\end{alignat*}
for $\CF \in \D_c(\FX)$ (Proposition \ref{propoext}). Both facts imply that  the ringed topos obstruction theory holds for locally noetherian formal schemes. In particular, from Corollary \ref{cordeformalsch} we get that the existence of a deformation $\FX$ of  $\FX'\xto{f'} \FY'$ by a closed immersion $\FY' \inc \FY$ given by a square zero ideal   $\CN$ is determined by an obstruction which lives in $ \ext^{2}_{\CO_{\FX'}}(\WCL_{\FX'/\FY'}, f'^{*} \CN)$. Uniqueness is controlled by $\ext^{1}_{\CO_{\FX'}}(\WCL_{\FX'/\FY'}, f'^{*} \CN)$. Moreover, obstruction conditions for the existence and uniqueness of lifting morphisms are given in Corollary \ref {cordemorfisformalsch}.  Explicitly, the existence of a deformation is guaranteed by the vanishing of an element in an $\ext^{1}$ group, and uniqueness is determined by an $\ext^0$ group.

Section \ref{sec5} is devoted to the characterizations of infinitesimal conditions through the cotangent complex. In Proposition \ref{cotanliso} (Corollary \ref{cotanetale}) we show that a morphism is smooth (\'etale) if, and only  if,  $\WCL_{\FX/\FY} \simeq \om^{1}_{\FX/\FY}[0]$ and $\om^{1}_{\FX/\FY}$ is a finite-rank locally free $\CO_{\FX}$-module (resp. $\WCL_{\FX/\FY} \simeq 0$).  A first consequence is that, for a class of non smooth usual schemes we can obtain deformations results depending on the differential module of a formal scheme. 
Corollary \ref{corclosedsmooth} provides that if  $j \colon \FX' \inc \FX$ is a closed immersion given by an  ideal $\CJ\subset \CO_{\FX}$ and $f \colon \FX \to \FY$ a smooth morphism, then 
$\tau^{\scriptscriptstyle{\ge -1}}(\WCL_{\FX'/\FY})  \simeq
(0 \to \CJ/\CJ^{2} \to j^*\om^{1}_{\FX/\FY}\to 0)$.
As an application  of these results, we obtain in Proposition  \ref{complinters} that the cotangent complex of a complete intersection morphism $\FX' \overset{j} \inc \FX \to \FY$ ($j$ a regular immersion and $\FX$ is $\FY$-smooth) is isomorphic to a complex of locally free $\CO_{\FX'}$-modules concentrated in degree $[-1,0]$:
\[
\WCL_{\FX'/\FY}  \simeq (0 \to \CJ/\CJ^{2} \to i^*\om^{1}_{\FX/\FY}\to 0).
\]

\begin{ack}

I would like  to thank  Leo Alonso and Ana Jerem\'{\i}as for their encouragement during the elaboration of this paper. Their comments about Grothendieck duality and Greenleess-May duality on formal schemes have illuminated my search until finding a suitable definition of cotangent complex.

I am also grateful to Luc Illusie for indicating me the lack of a treatment of cotangent complex for formal schemes and for his hospitality and support during my stay at Orsay, when it was my first contact with cotangent complex.
\end{ack}

\section{Preliminaries on formal schemes} \label{sec1}
We denote by $\sch$ the category of locally noetherian schemes and by $\sfn$ the category  of
locally noetherian formal schemes. The affine noetherian formal schemes are a full subcategory  of $\sfn$, denoted by $\sfna$. We  assume that the reader is familiarized with the theory of formal schemes as explained in \cite[\S 10]{EGA1}, \cite{AJL2}, \cite{AJP1} and \cite{AJP2}. However, in this section we fix some notations and recall some terminology that will be frequently used in this work.

\begin{parraf}\label{notlim}
Given a morphism $f:\FX \to \FY$ in $\sfn$ and $\CK \subset \CO_{\FY}$ an ideal of  definition, there exists an ideal of  definition $\CJ \subset \CO_{\FX}$ such that  $f^{*}(\CK) \CO_{\FX} \subset \CJ$ (see \cite[(10.5.4) and (10.6.10)]{EGA1}). 
For any such pair of ideals and for each $n \in \NN$, if  $f_{n}\colon X_{n}:=(\FX,\CO_{\FX}/\CJ^{n+1}) \to Y_{n}:=(\FY,\CO_{\FY}/\CK^{n+1})$ is the morphism in $\sch$ induced by $f$,  then $f$ is expressed as $f = \dirlim {n} f_{n}$.

\begin{subparraf}\label{existidefmorf}
The morphism $f$ is of \emph{pseudo finite type}  \cite[p. 7]{AJL2} (\emph{separated}  \cite[\S 10.15]{EGA1} and \cite[1.2.2]{AJL2}) if there exist $\CJ \subset \CO_{\FX}$ and  $\CK\subset \CO_{\FY}$ ideals of  definition  with  $f^{*}(\CK)\CO_{\FX} \subset \CJ$ (and hence for any pair of ideals satisfying this condition) such that $f_{0}$ is of  finite type (resp. separated). The morphism $f$ is of finite type if it is adic and of pseudo finite type (\cf\  \cite[(10.13.1)]{EGA1}).

Any (pseudo) finite type morphism $\spf(B) \to \spf(A)$ in $\sfna$ can be factorized as
\[
\spf(B) \overset{j} \inc  \spf(A\{ \mathbf{X}\}) \xto{p} \spf(A)
\]
\[
(\spf(B) \overset{j} \inc  \spf(A\{ \mathbf{X}\}[[\mathbf{Y}]] )\xto{p} \spf(A)),
\]
where $\mathbf{X}$ (and $\mathbf{Y}$) are finite numbers of indeterminates, $A\{ \mathbf{X}\}$ is the restricted formal power series ring, $j$ is a  closed immersion and $p$ is  the natural projection given by the morphism of adic rings $A \to A\{ \mathbf{X}\}$ (resp. $A \to A\{ \mathbf{X}\}[[\mathbf{Y}]]$). Any (pseudo) finite type morphism in $\sfn$ admits a local description as the one above \cite[10.13.1]{EGA1}(resp. \cf\ \cite[Proposition 1.7]{AJP1}).

Moreover, if $f \colon \spf(B) \to \spf(A)$ is a pseudo finite type morphism  in $\sfna$, then it can be expressed as the completion $f=\widehat{h}$ of a  finite type morphism $h \colon\spec(C) \to \spec(A)$ in $\scha$.
\end{subparraf}

\begin{subparraf}  \label{closed}
Given $\CJ \subset \CO_{\FX}$ an ideal and $\FX' \subset \FX$ the closed support of $\CO_{\FX}/\CJ$, $(\FX', (\CO_{\FX}/\CJ)|_{\FX'})$ is a formal scheme \cite[10.14.1, 10.14.2]{EGA1}. We say that $\FX'$ is the   \emph{closed formal subscheme} of  $\FX$ defined by $\CJ$.
A morphism $f:\FZ \to \FX$  is a \emph{closed immersion}  if there exists a closed formal subscheme  $\FX'\subset \FX$ such that $f$ factors as
$
\FZ \xto{g} \FX' \inc \FX
$
where $g$ is a isomorphism \cite[\S 10.14.2]{EGA1}. 
\end{subparraf}

\end{parraf}

\begin{parraf} \cite[Definition 2.1 and Definition 2.6]{AJP1} 
A morphism $f:\FX \to \FY$ in $\sfn$ is \emph{smooth  (unramified, \'etale)} if it is of pseudo finite type and satisfies the following lifting condition:

\emph{For all affine $\FY$-schemes $Z$ and for each closed subscheme $T\inc Z$ given by a square zero ideal $\CJ \subset \CO_{Z}$, the induced map
\begin{equation*} 
\Hom_{\FY}(Z,\FX) \lto \Hom_{\FY}(T,\FX)
\end{equation*}
is surjective (resp. injective, bijective)}.

If $f$ is of finite type, then we say that $f$ is \emph{adic smooth   (adic unramified, resp. adic \'etale)}.
\end{parraf}

\begin{parraf} \cite[(10.10.1)]{EGA1} \label{deftriangulito}
Let $\FX=\spf(A)$ with $A$ a  $I$-adic noetherian ring, $X=\spec(A)$ and $X'=\spec(A/I)$, so  $\FX= X_{/X'} \xto{\iota} X$ is a completion morphism. 
There is an additive covariant  functor from  the category of  $A$-modules to the category of  $\CO_{\FX}$-modules
\begin{equation} 
\begin{array}{ccc}
A\modu& \overset{\tr} \lto & \Modu(\FX)\\
M					&  \rightsquigarrow    &  M^{\tr}\\
\end{array}
\end{equation}
defined as
$
M^{\tr}:=(\widetilde{M})_{/X'} = 
\invlim {n} \frac{\widetilde{M}}{\widetilde{I}^{n+1}\widetilde{M}}$ and in the obvious way for morphisms. In the category of finite type $A$-modules it holds that $\iota^*(\widetilde{M})\simeq M^{\tr}$ and $\tr$ induces a category equivalence to coherent $\CO_{\FX}$-modules. 
\end{parraf}

\begin{parraf} 
Given $f: \FX \to \FY$ in $\sfn$ the  \emph{differential pair of  $\FX$ over $\FY$}, $( \om^{1}_{\FX/\FY}, \hd_{\FX/\FY})$, is 
locally  given  by
\[
\left( (\om^{1}_{B/A})^{\tr},  \CO_{\FU}=B^{\tr} \xto{\text{ via }\hd_{B/A}} (\om^{1}_{B/A})^{\tr} \right)
\] 
for all open subsets $\FU=\spf(B) \subset \FX$ and $\FV=\spf(A) \subset \FY$ with $f(\FU) \subset \FV$, being $\hd_{B/A}$ the completion of the canonical derivation $d_{B/A}$ w.r.t the adic topology in $B$. 
The   $\CO_{\FX}$-module $\om^{1}_{\FX/\FY}$ is called the \emph{module of  $1$-differentials of  $\FX$ over $\FY$} and the continuous $\FY$-derivation $\hd_{\FX/\FY}$ is called the \emph{canonical derivation of  $\FX$ over $\FY$}. Whenever $f$ is a pseudo finite type map, $\om^{1}_{\FX/\FY}$ is a coherent $\CO_{\FX}$-module. If $f$ is smooth (unramified), then $\om^{1}_{\FX/\FY}$ is a  locally free $\CO_{\FX}$-module of finite-rank (resp. $\om^{1}_{\FX/\FY}=0$).
The basic properties of  the differential pair in $\sfn$ are treated, for instance, in \cite[\S3, \S4]{AJP1} or \cite[\S2.6]{LNS}.

\end{parraf}


\section{The affine case} \label{sec2}

Let $A$ be an adic noetherian ring. We  denote by $\C(A)$ the category of complexes of $A$-modules and by $\D(A)$  its derived category. Let $\D^{-}(A)$ ($\D_{\cc}^{-}(A)$) be the full subcategory of $\D(A)$ whose complexes have bounded above cohomology (resp. coherent bounded above cohomology).  
Given $I \subset A$ an ideal of definition, the $I$-adic completion of modules induces the $I$-adic completion functor $\LA_{A, I} \colon \C(A) \to \C(A)$. Notice that $\LA_{A, I}$ only depends on the adic topology of $A$, not on the chosen ideal of definition, so once we  fix a  topology, we will denote it by $\LA_{A}$ or $\LA$.

Via flat resolutions $\LA$ has a left derived functor  $\LD \LA \colon \D^{-}(A) \to \D^{-}(A)$, the  \emph{derived completion functor}. Indeed, by \cite[p. 53, Theorem 5.1]{H} it suffices to prove that if $P \in  \D^{-}(A)$ is a flat acyclic complex, then $\LA(P)$ is acyclic.
By  \cite[(\textbf{0}.13.2.3)]{EGA31} one is reduced to showing that for all $n>0, P \otimes A/I^n$ is acyclic. But given $Q \xto{\sim} A/I^n$ a flat resolution, $P \otimes Q  \simeq P \otimes A/I^n $ and $P \otimes Q$ is acyclic, whence the conclusion.

\begin{parraf} \label{propicotang}
For a continuous map $A \to B$  of adic noetherian rings and $F \in \D^{-}(A)$, we summarize here some known properties that will be used on what follows:
\begin{enumerate}
\item \label{propicotang1}
$\LD \LA(\LD \LA(F)) = \LD \LA (F)$.

\item \label{propicotang3}
The canonical map
$
\LD \LA_{B}( F \underset{=}{\otimes} B) \to \LD \LA_B(\LD \LA_A( F) \underset{=}\otimes B)
$
is an isomorphism, where the symbol $\underset{=}{\otimes}$ stands for the derived tensor product.
\item \label{propicotang2}
 $F \in \D_{\cc}^{-} (A)$, then  $\LD \LA(F) \simeq F$.
\end{enumerate}
We give an outline of the proof: take a free resolution $P \xto{\sim} F$ in $\D^{-}(A)$. For (1) use that completion preserves flatness and is idempotent;  (2)  is a consequence of the module isomorphism $ \LA_{B}( F \otimes B) \simeq \LA_B(\LA_A( F) \otimes B)$; for (3) we can  take $P$ a resolution of  finite-rank free  modules.  Then apply the fact that  finite type modules are complete.
\end{parraf}

\begin{parraf}
Let $A \to B$  be a continuous map of adic noetherian rings. We denote by  $\WL_{B/A}:=\LA(L_{B/A}) \in \C^{\le 0}(B)$ the completion w.r.t the adic topology of $B$ of the usual cotangent complex $L_{B/A}$ \cite[Chapitre II, p. 123, 1.2.3]{I}.  We will call it   \emph{complete cotangent complex of $B$ over $A$}.

Observe that, since $L_{B/A}$ is a complex of flat $B$-modules \cite[Chapitre II, 1.2.3, 1.2.5]{I},  we can identify $\LD \LA(L_{B/A})$ with $\WL_{B/A}$ in $\D^{-}(B)$. Henceforth,  in $\D^{-}(B)$ we will use both notations for our convenience. 

The augmentation map in the usual case $L_{B/A} \to \Omega^{1}_{B/A}$ \cite[Chapitre II, 1.2.4]{I}, through the completion functor,  leads  to the augmentation map in $\C^{\le 0}(B)$
\begin{equation} \label{aug}
\WL_{B/A} \to \om^{1}_{B/A}.
\end{equation}
\end{parraf}

\begin{parraf}
As in the usual case \cite[Chapitre II, 1.2.3]{I}, the complete cotangent complex and the augmentation map depend functorially on the ring morphism. That is, given a commutative diagram of adic rings
\[
\begin{tikzpicture}
\draw[white] (0cm,2cm) -- +(0: \linewidth)
node (21) [black, pos = 0.43] {$A$}
node (22) [black, pos = 0.57] {$B$};
\draw[white] (0cm,0.5cm) -- +(0: \linewidth)
node (11) [black, pos = 0.43] {$A'$}
node (12) [black, pos = 0.57] {$B'$};
\draw [->] (21)  --   (22);
\draw [->] (11)  --   (12);
\draw [->] (21)  --   (11);
\draw [->] (22)  --   (12);
\end{tikzpicture}
\]
there is a morphism  $\WL_{B/A} \to \WL_{B'/A'}$ in $\C^{\le 0}(B)$ such that the following diagram is commutative:
\[
\begin{tikzpicture}
\draw[white] (0cm,2cm) -- +(0: \linewidth)
node (21) [black, pos = 0.4] {$\WL_{B/A}$}
node (22) [black, pos = 0.6] {$\WL_{B'/A'}$};
\draw[white] (0cm,0.5cm) -- +(0: \linewidth)
node (11) [black, pos = 0.4] {$\om^{1}_{B/A}[0]$}
node (12) [black, pos = 0.6] {$\om^{1}_{B'/A'}[0]$.};
\draw [->] (21)  --   (22);
\draw [->] (11)  --   (12);
\draw [->] (21)  --   (11);
\draw [->] (22)  --   (12);
\end{tikzpicture}
\]
Note that the morphism $\WL_{B/A} \to \WL_{B'/A'}$ is equivalent (via adjunction) to giving a morphism $\WL_{B/A}\otimes_{B}B' \to \WL_{B'/A'}$, or for that matter $L_{B/A}\tc_{B}B' \to \WL_{B'/A'}$ in $\C^{\le 0}(B')$.
\end{parraf}
\begin{parraf} \label{llanpropis}
Let  $A \to B $ and $B \to C $ be two morphisms of adic noetherian rings. Since $\LD \LA$ is a $\Delta$-functor, using \ref{propicotang}.(\ref{propicotang3})  the associated distinguished triangle \cite[Chapitre II, Proposition 2.1.2]{I} leads to a distinguished triangle  in $\D^{-}(C)$ 
\[
\LD \LA (\WL_{B/A} \otimes
_{B}C ) \to \WL_{C /A} \to \WL_{C /B} \overset{+} \to 
\]
such that the maps in degree zero are given by the natural morphisms $L_{B/A} \tc_{B}C  \to \WL_{C /A} \to \WL_{C /B}$ in $\C^{\le 0}(C)$.
\end{parraf}

The next proposition is a re-statement of \cite[Proposition 7.1.29]{GR} in the $\sfna$ context, and the proof we give is a rewriting of the proof of the above mentioned result. 
\begin{propo} \label{smoadic}
If $f \colon \spf(B) \to \spf(A)$  is an \emph{adic} smooth morphism in $\sfna$, then the augmentation map (\ref{aug}) gives an isomorphism $\WL_{B/A} \simeq
\om^{1}_{B/A}[0]$ and $\WL_{B/A} \in \D^{-}_{\cc}(B)$. 
\end{propo}

\begin{proof}
Let $I \subset A$ be an ideal of definition. For each $n \in \NN$, put $A_{n}:= A/I^{n+1}$ and $B_{n}:= B/(IB)^{n+1}$. 
Since $f$ is flat \cite[Proposition 4.8]{AJP1}, by base-change \cite[Chapitre II, Corollaire 2.2.3]{I} we have that  $L_{B/A}\otimes_{B} B_{n} \simeq
L_{B_{n}/A_{n}},\,  \forall n$. 
By  base-change \cite[Corollary 5.2]{AJP2}, each $f_{n}: \spec (B_{n}) \to \spec (A_{n})$ is a smooth morphism in $\sch$, therefore the augmentation map $L_{B_{n}/A_{n}} \to \Omega^{1}_{B_{n}/A_{n}}$ gives an isomorphism $L_{B_{n}/A_{n}} \simeq
\Omega^{1}_{B_{n}/A_{n}}[0]$ (\cite[Chapitre III, Proposition 3.1.2]{I}, \cite[Chapitre II, 2.3.6.3]{I}).

From \cite[\textbf{0}.13.2.3]{EGA31} it follows that $\h^{i}(\WL_{B/A})=0$, for all $i \le -1$,
and $\h^{0}(\WL_{B/A}) \simeq \invlim  {n} \h^{0}(L_{B_{n}/A_{n}})=\invlim  {n} \Omega^{1}_{B_{n}/A_{n}} = \om^{1}_{B/A}$,  by \cite[\textbf{0}.20.7.14]{EGA41}.

The last assertion follows from the coherence of the differential module \cite[Proposition 3.3]{AJP1}.
\end{proof}

\begin{cor} \label{etadic}
If $f \colon \spf(B) \to \spf(A)$  is an \emph{adic} \'etale morphism in $\sfna$, then $\WL_{B/A} \simeq
0$.
\end{cor}

\begin{proof}
It's immediate from the last proposition and  \cite[Proposition 4.6]{AJP1}.
\end{proof}

\begin{propo} \label{cloim}
Let $j \colon \spf(C) \inc \spf(B)$ be a closed immersion in $\sfna$. Then $\WL_{C/B} \simeq  L_{C/B}$. In particular:
\begin{enumerate}
\item
$\WL_{C/B} \in \D^{-}_{\cc}(C)$.
\item
$\h^{0}(\WL_{C/B})=0$.
\item
If $J \subset B$ is the ideal defining $j$, then there exists a canonical functorial isomorphism $\h^{1}(\WL_{C/B}) \xto{\sim}  J/J^{2}$.
\end{enumerate}
\end{propo}

\begin{proof}
Since the morphism $B \to C$ is surjective,  we have that $L_{C/B} \in \D^{-}_{\cc}(C)$ \cite[Chapitre II, Proposition 2.3.6.3, Corollaire 2.3.7]{I}, then by \ref{propicotang}.(\ref{propicotang2})  $\WL_{C/B} \simeq L_{C/B}$. Then (1)  is immediate and the other two assertions  follow from \cite[Chapitre III, Corollaire 1.2.8.1]{I}.
\end{proof}

\begin{lem} \label{compl}
Let $A$ be a noetherian ring, $I \subset A$ an ideal  and $\HA$ the $I$-adic completion ring. Then $\WL_{\HA /A}\simeq
0$.                                                                                                                                                                                                                                                                                                                                                                                                                                                                                                                                                                                                                                                                                                                                                                                                                                                                                                                                                                                                                                                  
\end{lem}

\begin{proof}
Since the canonical  map $A \to \HA$ is flat,  base-change  \cite[Chapitre II
, p. 138, Proposition 2.2.1]{I} implies that $L_{\HA/A} \otimes_{\HA} \HA/I^{n}\HA \simeq
0, \forall n>0$.
The result is immediate from \cite[\textbf{0}.13.2.3]{EGA31}.
\end{proof}

\begin{lem} \label{espafin}
Let $A$ be an adic noetherian ring, $\mathbf{X}$ and $\mathbf{Y}$ finite numbers of indeterminates and $A \to A\{ \mathbf{X}\}[[\mathbf{Y}]]$ the canonical morphism. Then the augmentation map (\ref{aug}) gives an isomorphism
$ \WL_{A\{ \mathbf{X}\}[[\mathbf{Y}]]/A} \simeq
\om^{1}_{A\{ \mathbf{X}\}[[\mathbf{Y}]]/A}[0]$ and $ \WL_{A\{ \mathbf{X}\}[[\mathbf{Y}]]/A} \in \D^{-}_{\cc}(A\{ \mathbf{X}\}[[\mathbf{Y}]])$.
\end{lem}

\begin{proof}
The canonical morphism $A \to A\{ \mathbf{X}\}[[\mathbf{Y}]]$ admits a factorization 
\[
A \to A\{\mathbf{X, Y}\}  \to A\{\mathbf{X}\}[[\mathbf{Y}]].
\]
By \ref{llanpropis} we get a distinguished triangle 
\[
\LD \LA(\WL_{A\{ \mathbf{X,Y}\}/A} \otimes
_{A\{ \mathbf{X,Y}\}}A\{ \mathbf{X}\}[[\mathbf{Y}]] )\to \WL_{A\{ \mathbf{X}\}[[\mathbf{Y}]]/A} \to \WL_{A\{ \mathbf{X}\}[[\mathbf{Y}]]/A\{ \mathbf{X,Y}\}} \overset{+} \to
\]
in $\D^{-}(A\{ \mathbf{X}\}[[\mathbf{Y}]])$.
It  follows that
\begin{alignat*}{3}
\WL_{A\{ \mathbf{X}\}[[\mathbf{Y}]]/A}&\simeq
\LD \LA(\WL_{A\{ \mathbf{X,Y}\}/A} \otimes
_{A\{ \mathbf{X,Y}\}}A\{ \mathbf{X}\}[[\mathbf{Y}]] )&  \text{Lemma \ref{compl}}\\    
&\simeq
 \WL_{A\{ \mathbf{X,Y}\}/A} \otimes_{A\{ \mathbf{X,Y}\}}A\{ \mathbf{X}\}[[\mathbf{Y}]] 	&  \text{Prop. \ref{smoadic}, \ref{propicotang}.(\ref{propicotang2})}\\
&\simeq
\om^{1}_{A\{ \mathbf{X,Y}\}/A}[0] \otimes_{A\{ \mathbf{X,Y}\}}A\{ \mathbf{X}\}[[\mathbf{Y}]]&\quad \text{Prop. \ref{smoadic}} \\
& \simeq \om^{1}_{A\{ \mathbf{X}\}[[Y]]/A}[0]&\quad \text{\cite[Cor. 4.10]{AJP1}}
\end{alignat*}
By coherence of the differential module \cite[Proposition 3.3]{AJP1} we obtain that $\WL_{A\{ \mathbf{X}\}[[\mathbf{Y}]]/A} \in \D^{-}_{\cc}(A\{ \mathbf{X}\}[[\mathbf{Y}]])$.
\end{proof}

\begin{propo} \label{pftcoh}
If $f \colon \spf(B) \to \spf(A)$  is a pseudo  finite type morphism in $\sfna$, then $\WL_{B/A} \in \D^{-}_{\cc}(B)$.
\end{propo}

\begin{proof}
The morphism $A \to B$ admits a factorization 
\begin{equation}\label{fact1}  
A \to A\{ \mathbf{X}\}[[\mathbf{Y}]]  \to B \to 0
\end{equation}
where $\mathbf{X}$ and $\mathbf{Y}$ are finite numbers of indeterminates (see \ref{existidefmorf}), from which we obtain a distinguished triangle in $\D^{-}(B)$
\[
\LD \LA(\WL_{A\{ \mathbf{X}\}[[\mathbf{Y}]]/A} \otimes 
_{A\{ \mathbf{X}\}[[\mathbf{Y}]]}B) \to \WL_{B/A} \to \WL_{B/A\{ \mathbf{X}\}[[\mathbf{Y}]]} \overset{+} \to.
\]
Lemma \ref{espafin}, \ref{propicotang}.(\ref{propicotang2}) and  Proposition \ref{cloim} show that the left and the right terms in the triangle are in $\D^{-}_{\cc}(B)$, so the middle term is also in $\D^{-}_{\cc}(B)$. 
\end{proof}

\begin{cor} \label{trian}
Let $f \colon \spf(C) \to \spf(B)$ and $g \colon\spf(B) \to \spf(A)$ be two pseudo finite morphisms in $\sfna$. Then there exists a distinguished triangle in $\D^{-}_{\cc}(C)$
\begin{equation} \label{distrian}
\WL_{B/A} \otimes_{B}C  \to \WL_{C /A} \to \WL_{C /B} \overset{+} \to
\end{equation}
such that the  maps  in degree zero are induced by the canonical morphisms in $\C^{\le 0}(C)$.

We will call the distinguished triangle (\ref{distrian}) the \emph{complete distinguished triangle} associated to $g \circ f \colon\spf(C) \to \spf(B)  \to \spf(A)$.
\end{cor}

\begin{proof}
By  Proposition \ref{pftcoh} and \ref{propicotang}.(\ref{propicotang2})  the triangle of  \ref{llanpropis} is rewritten in this form.
\end{proof}

\begin{cor} \label{loc}
Let $f \colon \spf(C) \to \spf(B)$  and $g \colon \spf(B) \to \spf(A)$ two pseudo  finite type morphisms in $\sfna$. 
If $f$ is an adic \'etale or a completion morphism, then  the  morphism $\WL_{B/A} \otimes_{B} C \to  \WL_{C/A}$ is an isomorphism.
\end{cor}

\begin{proof}
Corollary \ref{etadic} and Lemma \ref{compl} imply that $\WL_{C/B} \simeq 0$. Hence the result follows from the associated complete distinguished triangle.
\end{proof}

The next proposition upgrades Proposition \ref{smoadic}, and uses in a critical way the results for completions in Corollary \ref{loc}.

\begin{propo} \label{smo}
If $f \colon \spf(B) \to \spf(A)$  is a  smooth morphism in $\sfna$, then the augmentation map (\ref{aug}) gives an isomorphism $\WL_{B/A} \simeq
\om^{1}_{B/A}[0]$.
\end{propo}

\begin{proof}
By Corollary \ref{loc}, \cite[Proposition 4.1 and 4.3]{AJP1} and \cite[Theorem 7.13]{AJP2} we may assume that $f$ factorizes as 
\begin{equation} \label{eq}
\spf(B) \xto{\kappa} \spf(B')\xto{f'} \spf(A)
\end{equation}
where $\kappa$ is a morphism of completion (that is, $B$ is an adic completion of $B'$) and $f'$ is  adic smooth. 
The result follows from the canonical commutative diagram in $\D^{-}_{\cc}(B)$ 
\[
\begin{tikzpicture}
\draw[white] (0cm,2cm) -- +(0: \linewidth)
node (21) [black, pos = 0.35] {$ \WL_{B'/A} \otimes_{B'}B$}
node (22) [black, pos = 0.65] {$ \WL_{B /A}$};
\draw[white] (0cm,0.5cm) -- +(0: \linewidth)
node (11) [black, pos = 0.35] {$(\om^{1}_{B'/A} \otimes_{B'}B)[0]$}
node (12) [black, pos = 0.65] {$\om ^{1}_{B /A}[0] $.};
\draw [->] (21)  --   (22) node[above, midway, scale=0.75]{$\simeq$};
\draw [->] (11)  --   (12)node[above, midway, scale=0.75]{$\simeq$};
\draw [->] (21)  --   (11)node[left, midway, scale=0.75]{$\simeq$};
\draw [->] (22)  --   (12);
\end{tikzpicture}
\]
Being $f'$ adic and smooth, the left vertical map is an isomorphism by Proposition \ref{smoadic}. Since $\kappa$ is a completion morphism,   the horizontal arrows are isomorphisms, by Corollary \ref{loc} and \cite[Corollary 4.10]{AJP1}.
\end{proof}

\begin{cor} \label{et}
If $f \colon \spf(B) \to \spf(A)$  is an \'etale morphism in $\sfna$, then $\WL_{B/A} \simeq
0$.
\end{cor}

\begin{proof}
It is immediate applying the last proposition and \cite[Proposition 4.6]{AJP1}.
\end{proof}

\begin{parraf} \label{cotsef2}
Let $j \colon \spf(C) \inc \spf(B)$ be a closed immersion given by an ideal $J\subset B$ and $f \colon \spf(B)  \to \spf(A)$ a pseudo finite type morphism in $\sfna$. Then there is a commutative diagram of coherent modules
\begin{equation}\label{comdiagh1h0}
\begin{tikzpicture}
\draw[white] (0cm,2cm) -- +(0: \linewidth)
node (21) [black, pos = 0.35] {$ \h^{1}(\WL_{C/B})$}
node (22) [black, pos = 0.65] {$ \h^{0}(\WL_{B/A} \otimes_{B}C)$};
\draw[white] (0cm,0.5cm) -- +(0: \linewidth)
node (11) [black, pos = 0.35] {$J/J^2$}
node (12) [black, pos = 0.65] {$\om^{1}_{B/A}\otimes_{B}C$};
\draw [->] (21)  --   (22);
\draw [->] (11)  --   (12)node[above, midway, scale=0.75]{$-\hd$};
\draw [->] (21)  --   (11)node[left, midway, scale=0.75]{$\simeq$} node[right, midway, scale=0.75]{\ref{cloim}};
\draw [->] (22)  --   (12)node[right, midway, scale=0.75]{(\ref{aug})};
\end{tikzpicture}
\end{equation}
where the top  row is  the connection morphism of the cohomology sequence given by the complete distinguished triangle associated to $f \circ j$ and $\hd$ is induced by  the complete canonical derivation $\hd_{B/A}$.
Indeed, there is  a commutative diagram in $\D^{-}(B)$:

\[
\begin{tikzpicture}
\draw[white] (0cm,2cm) -- +(0: \linewidth)
node (21) [black, pos = 0.2] {$  L_{B/A} \otimes_{B}C$}
node (22) [black, pos = 0.4] {$ L_{C/A}$}
node (23) [black, pos = 0.6] {$ L_{C/B}$}
node (24) [black, pos = 0.8] {};
\draw[white] (0cm,0.5cm) -- +(0: \linewidth)
node (11) [black, pos = 0.2] {$\WL_{B/A} \otimes_{B}C$}
node (12) [black, pos = 0.4] {$\WL_{C/A}$}
node (13) [black, pos = 0.6] {$\WL_{C/B}  $}
node (14) [black, pos = 0.8] {};
\draw [->] (21)  --   (22);
\draw [->] (22)  --   (23);
\draw [->] (23)  --   (24)node[above, midway, scale=0.75]{$\delta$} node[below, midway, scale=0.75]{$+$};
\draw [->] (11)  --   (12);
\draw [->] (12)  --   (13);
\draw [->] (13)  --   (14)node[above, midway, scale=0.75]{$\LD \LA (\delta)$} node[below, midway, scale=0.75]{$+$};
\draw [->] (21)  --   (11);
\draw [->] (22)  --   (12);
\draw [->] (23)  --   (13)node[left, midway, scale=0.75]{$\simeq$} node[right, midway, scale=0.75]{\ref{cloim}};
\end{tikzpicture}
\]
where the top row is the usual distinguished triangle \cite[Chapitre II, Proposition 2.1.2]{I} and the bottom row is the complete distinguished triangle (\ref{distrian}).
 By \cite[Chapitre III, Proposition 1.2.9]{I}  we get the  diagram 
\[
\begin{tikzpicture}
\draw[white] (0cm,3.5cm) -- +(0: \linewidth)
node (41) [black, pos = 0.2] {$ \h^{1}(L_{C/B})$}
node (43) [black, pos = 0.6] {$ \h^{0}(L_{B/A} \otimes_{B}C)$};
\draw[white] (0cm,2.5cm) -- +(0: \linewidth)
node (32) [black, pos = 0.4] {$ \h^{1}(\WL_{C/B})$}
node (33) [black, pos = 0.6] {}
node (34) [black, pos = 0.8] {$ \h^{0}(\WL_{B/A} \otimes_{B}C)$};
\draw[white] (0cm,1.5cm) -- +(0: \linewidth)
node (21) [black, pos = 0.2] {$ J/J^2$}
node (22) [black, pos = 0.4] {}
node (23) [black, pos = 0.6] {$\Omega^{1}_{B/A}\otimes_{B}C$};
\draw[white] (0cm,0.5cm) -- +(0: \linewidth)
node (12) [black, pos = 0.4] {$J/J^2$}
node (14) [black, pos = 0.8] {$\om^{1}_{B/A}\otimes_{B}C$};
\draw [->] (41)  --   (43)node[above, midway, scale=0.75]{via $\delta$};
\draw [->] (32)  --   (34)node[above, near start, scale=0.75]{via $\LD \LA (\delta)$};
\draw [-] (21)  --   (22)node[above, near end, scale=0.75]{};
\draw [->] (22)  --   (23)node[above, midway, scale=0.75]{$-d$};
\draw [->] (12)  --   (14)node[above, midway, scale=0.75]{$-\hd$};
\draw [->] (41)  --   (21)node[left, midway, scale=0.75]{$\simeq$};
\draw [->] (32)  --   (12)node[left, near start, scale=0.75]{$\simeq$};
\draw [-] (43)  --   (33) node[right, midway, scale=0.75]{(\ref{aug})};
\draw [->] (33)  --   (23)node[scale=0.75]{};
\draw [->] (34)  --   (14)node[right, midway, scale=0.75]{(\ref{aug})};
\draw [->] (41)  --   (32)node[above, midway, scale=0.75]{$\simeq$};
\draw [->] (21)  --   (12)node[above, midway, scale=0.75]{$=$};
\draw [->] (43)  --   (34);
\draw [->] (23)  --   (14);
\end{tikzpicture} 
\]
where,  as a consequence of the commutativity of the   others squares in the diagram, the square in front, namely (\ref{comdiagh1h0}), is also commutative.
\end{parraf}

\begin{propo} \label{h0afin}
If $f \colon \spf(B) \to \spf(A)$ is a pseudo finite type morphism in $\sfna$, then the augmentation map (\ref{aug}) induces an isomorphism 
\[
\h^{0}(\WL_{B/A}) \simeq\om^{1}_{B/A}.
\]
\end{propo}

\begin{proof}
The morphism $A \to B$ admits a factorization (see \ref{existidefmorf})
\begin{equation} \label{factoriz}
A \to A\{ \mathbf{X}\}[[\mathbf{Y}]]  \to B \to 0
\end{equation}
where $\mathbf{X}$ and $\mathbf{Y}$ denote two finite numbers of indeterminates. From the complete distinguished triangle, 
by   \ref{cotsef2}, we obtain a commutative diagram of $B$-modules
\[
\begin{tikzpicture}
\draw[white] (0cm,2cm) -- +(0: \linewidth)
node (21) [black, pos = 0.15] {$ \h^{1}(\WL_{B/A\{\mathbf{X}\}[[Y]]})$}
node (22) [black, pos = 0.5] {$ \h^{0}(\WL_{A\{\mathbf{X}\}[[Y]]/A} \otimes_{A\{\mathbf{X}\}[[Y]]}B)$}
node (23) [black, pos = 0.8] {$\h^{0}(\WL_{B /A})$}
node (24) [black, pos = 0.95] {$0$};
\draw[white] (0cm,0.5cm) -- +(0: \linewidth)
node (11) [black, pos = 0.15] {$J/J^2$}
node (12) [black, pos = 0.5] {$\om^{1}_{A\{\mathbf{X}\}[[Y]]/A} \otimes_{A\{\mathbf{X}\}[[Y]]}B$}
node (13) [black, pos = 0.8] {$\om^{1}_{B /A}$}
node (14) [black, pos = 0.95] {$0$};
\draw [->] (21)  --   (22);
\draw [->] (22)  --   (23);
\draw [->] (23)  --   (24);
\draw [->] (11)  --   (12);
\draw [->] (12)  --   (13);
\draw [->] (13)  --   (14);
\draw [->] (21)  --   (11) ;
\draw [->] (22)  --   (12);
\draw [->] (23)  --   (13);
\end{tikzpicture}
\]
where the top row is the end of the long exact sequence of cohomology given by the complete distinguished triangle associated to   (\ref{factoriz}) and the bottom row corresponds to the \emph{complete} second fundamental exact sequence (\cf\ \cite[Proposition 3.13]{AJP1} and \cite[\textbf{0}.20.7.17, \textbf{0}.20.7.20]{EGA41}).
Now, by Proposition \ref{cloim} and Lemma \ref{espafin}  we have that the two left maps  are isomorphisms, therefore so is  the right one.
\end{proof}

\section{Definition and properties} \label{sec3}

Given $(\FX, \CO_\FX)$ in $\sfn$,  denote by $\A(\FX)$  
the category of  $\CO_\FX$-modules and by $\D(\FX)$ its  derived category. Write $\A_\cc(\FX) \subset \A(\FX)$  ($\A_{\ic}(\FX)\subset \A(\FX)$ \cite[\S 3]{AJL2}) for the subcategory of coherent $\CO_\FX$-modules (resp. direct limits of coherent $\CO_{\FX}$-modules) and $\D_{\cc}(\FX)$ (resp. $\D_{{\ic}}(\FX)$) for the corresponding derived category. Let $\D^{-}(\FX)$ ($\D_{\cc}^{-}(\FX)$, $\D_{{\ic}}^{-}(\FX)$) be the full subcategory of $\D(\FX)$ whose complexes have bounded above cohomology (resp. coherent bounded above cohomology, coherent direct limit  bounded above cohomology).  
We write $\R$($\LD$) for the right(left)-derived  functors, that will be obtained via q-injective (resp. q-flat) resolutions.

\begin{parraf} 
\cite[1.2.1, 1.2.4, \S 5, \S 6]{AJL2}
Given $\FX \in \sfn$ and $\CI$ an ideal of definition of $\FX$, let  $\varGamma'_\FX \colon \A(\FX) \to \A(\FX)$ be the functor defined by
$\varGamma'_\FX \CF := \dirlim{n > 0} \shom_{\CO_\FX}(\CO_\FX/\CI^n, \CF)$. It does not
 depend on the  ideal of definition. We will denote by $\RGA_\FX$ the right derived functor of $\varGamma'_\FX $.

The homology localization functor
$\RLA_{\FX}\colon \D(\FX) \to \D(\FX)$,  defined as 
\[
\RLA_{\FX}(-):=\rshom^\bullet(\RGA_{\FX} \CO_{\FX}, -)
\]
is right adjoint to $\RGA_{\FX}$.
Whenever $\FX$ is   separated and $\CF $  is $q$-flat complex of $ \A_{\ic}(\FX)$-modules, 
\emph{Greenlees-May duality} \cite[Theorem 0.1]{AJL3} establishes  a functorial isomorphism 
\[
\RLA_{\FX}(\CF)\cong  \LA_{\FX}(\CF):=\invlim {n}(\CF\otimes\CO_{\FX}/\CI^n ).
\]
\end{parraf}

Below we summarize some known properties about  $\RGA$ and $\RLA$   that will be used along this work.  
 
\begin{parraf} \label{propihomologyloc}

Let $f \colon \FX \to \FY$ be a morphism in $\sfn$ and $\CE, \CF \in \D(\FX), \CG \in \D(\FY)$.
\begin{enumerate}
\item \label{gaf*ga}
The natural map $  \RGA_\FX \LD f^* \CG \to \RGA_\FX \LD f^* \RGA_\FY  \CG$ is an isomorphism (\cite[Proposition 5.2.8.(c)]{AJL2}).
\item \label{propihomlocmorfismoad}
If $f$ is  adic, then $\LD f^*\RGA_{\FY} \CG \cong \RGA_{\FX} \LD f^* \CG$ (\cite[Corollary 5.2.11.(c)]{AJL2}).
\item \label{homlocga}
$\RGA_{\FX} \RGA_{\FX} \simeq \RGA_{\FX}, \RLA_{\FX} \RLA_{\FX} \cong \RLA_{\FX}, \RGA_{\FX} \RLA_{\FX} \simeq \RGA_{\FX}, \RLA_{\FX} \RGA_{\FX} \cong \RLA_{\FX}  $ (\cite[Remarks 6.3.1.(1)]{AJL2}).
\item \label{homlala}
$\rshom^\bullet(\RLA_{\FX} \CE, \RLA_{\FX}\CF) \cong \rshom^\bullet( \CE,\RLA_{\FX} \CF) $ (\cite[Remarks 6.3.1.(1)]{AJL2}).
\item \label{galacoh}
If $ \CF$ in $\D_c(\FX)$, then $\rshom^\bullet( \CE, \CF) \cong\rshom^\bullet(\RGA_{\FX} \CE,\CF) $.
Therefore  $\RLA_\FX \CF\cong \CF$ (\cite[Proposition 6.2.1]{AJL2}).
\end{enumerate} 

\end{parraf}

\begin{lem} \label{lemderinver}
Given a morphism $f \colon \FX \to \FY$ in $\sfn$ and $\CG \in \D(\FY)$, there is a natural  isomorphism:
\[
\RLA_\FX \LD f^* \RLA_\FY \CG \xto{\sim} \RLA_\FX \LD f^* \CG
\]
\end{lem}

\begin{proof}
The morphism is the composition of the canonical isomorphisms: 
\begin{alignat*}{3}
\RLA_\FX \LD f^* \RLA_\FY \CE&\cong \RLA_\FX \RGA_\FX \LD f^*\RLA_\FY \CG&\text{\ref{propihomologyloc}.(\ref{homlocga})} \\
&\cong \RLA_\FX \RGA_\FX \LD f^* \RGA_\FY \RLA_\FY \CG& \quad\text{\ref{propihomologyloc}.(\ref{gaf*ga})}\\
&\cong\RLA_\FX \RGA_\FX \LD f^* \RGA_\FY \CG&\text{\ref{propihomologyloc}.(\ref{homlocga})}  \\
&\cong \RLA_{\FX} \RGA_\FX \LD f^{*} \CG&\text{\ref{propihomologyloc}.(\ref{gaf*ga})} \\
&\cong\RLA_{\FX} \LD f^{*} \CG &\text{\ref{propihomologyloc}.(\ref{homlocga})}
\end{alignat*}
\end{proof}

\begin{defn}
Let $f \colon \FX \to \FY$ be a morphism in $\sfn$. 
We define the \emph{complete cotangent complex of $\FX$ over $\FY$} or \emph{the complete cotangent complex of $f$} as $\WCL_{\FX/\FY}:= \RLA_{\FX}(\CL_{\FX/\FY}) \in \D^-(\FX)$, where $\CL_{\FX/\FY}$ is the usual cotangent complex of ringed spaces \cite[Chapitre II, 1.2.7]{I}.

Applying $\RLA_{\FX}$ to the composition of 
the  usual augmentation \cite[Chapitre II, (1.2.7.3)]{I} and the completion map $\Omega^1_{\FX/\FY} \to \om^1_{\FX/\FY}$, we obtain an augmentation map
$\WCL_{\FX/\FY} \to \RLA_{\FX}\om^1_{\FX/\FY} $ in $\D(\FX)$. Notice that, if $f$ is of pseudo finite type, $\om^1_{\FX/\FY} \in \A_\cc(\FX)$, and by \ref{propihomologyloc}.(\ref{galacoh}) $\RLA_{\FX}(\om^1_{\FX/\FY}) \cong\om^1_{\FX/\FY}$. Then for a pseudo finite type morphism the augmentation map can be written as
\begin{equation} \label{sheafaug}
\WCL_{\FX/\FY} \to  \om^1_{\FX/\FY}.
\end{equation}

\end{defn}

\begin{rem}
Since $\CL_{\FX/\FY}$ is isomorphic to a bounded above complex of $\CO_{\FX}$-free modules \cite[Chapitre II, (1.2.3.3)]{I},  it is a  $q$-flat complex of $ \A_{\ic}(\FX)$-modules. Therefore, by Greenlees-May duality, for $\FX$  separated it holds that $\WCL_{\FX/\FY}\cong\LA_{\FX}(\CL_{\FX/\FY})$.
\end{rem}

\begin{parraf}  \label{functoriality}
Given a commutative diagram of pseudo finite type morphisms in $\sfn$
\[
\begin{tikzpicture}
\draw[white] (0cm,2cm) -- +(0: \linewidth)
node (21) [black, pos = 0.43] {$\FX'$}
node (22) [black, pos = 0.57] {$\FY'$};
\draw[white] (0cm,0.5cm) -- +(0: \linewidth)
node (11) [black, pos = 0.43] {$\FX$}
node (12) [black, pos = 0.57] {$\FY$,};
\draw [->] (21)  --   (22)node[above, midway, scale=0.75]{$f'$};
\draw [->] (11)  --   (12)node[above, midway, scale=0.75]{$f$};
\draw [->] (21)  --   (11)node[left, midway, scale=0.75]{$u$};
\draw [->] (22)  --   (12);
\end{tikzpicture}
\] 
 applying  $\RLA_{\FX'}$ to the canonical map  $u^{*} \CL_{\FX/\FY} \to \CL_{\FX'/\FY'} $ (\cf\ \cite[Chapitre II, (1.2.7.2)]{I}), by Lemma \ref{lemderinver} we obtain a map\footnote{Notice that $\LD u^*=u^*$ for $q$-flat complexes.}  $\RLA_{\FX'} \LD u^{*} \WCL_{\FX/\FY} \to \WCL_{\FX'/\FY'}$ in $\D^-(\FX')$ such that the following diagram is commutative:
\[
\begin{tikzpicture}
\draw[white] (0cm,3.5cm) -- +(0: \linewidth)
node (21) [black, pos = 0.4] {$\RLA_{\FX'} \LD u^{*}\WCL_{\FX/\FY}$}
node (22) [black, pos = 0.6] {$\WCL_{\FX'/\FY'}$};
\draw[white] (0cm,2cm) -- +(0: \linewidth)
node (11) [black, pos = 0.4] {$\RLA_{\FX'}\LD u^{*} \om^{1}_{\FX/\FY} $};
\draw[white] (0cm,0.5cm) -- +(0: \linewidth)
node (01) [black, pos = 0.4] {$u^{*} \om^{1}_{\FX/\FY} $}
node (02) [black, pos = 0.6] {$\om^{1}_{\FX'/\FY'}$.};
\draw [->] (21)  --   (22);
\draw [->] (01)  --   (02);
\draw [->] (21)  --   (11);
\draw [->] (11)  --   (01);
\draw [->] (22)  --   (02);
\end{tikzpicture}
\]
 
\end{parraf}

\begin{rem}
We will show in Proposition \ref{coh} that if $f$  is pseudo finite type, then $ \WCL_{\FX/\FY} \in \D_{\cc}(\FX)$ and therefore $\RLA_{\FX'} \LD u^{*} \WCL_{\FX/\FY}\cong\LD u^{*} \WCL_{\FY/\FS}$. 
\end{rem}

\begin{propo} \label{triannocoh}
Let $f \colon \FX \to \FY$ and $g \colon \FY \to \FZ$ be two  morphisms in $\sfn$. Then there is a distinguished triangle in $\D^{-}(\FX)$
\[
\RLA_\FX \LD f^{*} \WCL_{\FY/\FS} \to \WCL_{\FX/\FS} \to \WCL_{\FX/\FY} \xto{+} 
\]
such that the maps in degree zero are the canonical ones.
\end{propo}

\begin{proof}
Since $\RLA_\FX$ is a $\Delta$-functor  it's sufficient to apply it  to the distinguished triangle $ f^{*} \CL_{\FY/\FS} \to \CL_{\FX/\FS} \to \CL_{\FX/\FY} \xto{+}  $ \cite[Chapitre II, 2.1.5.6]{I}.
\end{proof}

\begin{propo}\label{localcotag}
If $f \colon \FX \to \FY$ is a  morphism  in $\sfn$ and $i\colon \FU \inc \FX$, $j\colon \FV \inc \FY$  are two open immersions in $\sfn$ such that $f(\FU) \subset \FV$, then
\[
i^* \WCL_{\FX/\FY} \simeq \WCL_{\FU/\FV} \in \D^{-}(\FU).
\]
\end{propo}

\begin{proof}
Given a q-injective resolution $\CL_{\FX/\FY} \xto{\sim} \CI$   in $\D(\FX)$, $i^*\CL_{\FX/\FY} \xto{\sim} i^* \CI$ is also a q-injective resolution \cite[Lemma 2.4.5.2]{L}. It holds that: 
\begin{alignat*}{3}
i^* \WCL_{\FX/\FY}&=  i^*\rshom^\bullet(\RGA_{\FX} \CO_{\FX}, \CL_{\FX/\FY}) &  \\    
&=  i^*\shom^\bullet(\RGA_{\FX} \CO_{\FX}, \CI)& \\    
& =\shom^\bullet(i^*\RGA_{\FX} \CO_{\FX}, i^*\CI)& \\
& \simeq\shom^\bullet(\RGA_{\FU} \CO_{\FU}, i^*\CI)& \qquad \ref{propihomologyloc}.(\ref{propihomlocmorfismoad})\\
&=\rshom^\bullet(\RGA_{\FU} \CO_{\FU}, i^*\CL_{\FX/\FY})&\\
&=\rshom^\bullet(\RGA_{\FU} \CO_{\FU}, \CL_{\FU/\FV}) &\qquad \text{\cite[Chapitre II,(1.2.3.5)]{I}}\\
&=  \WCL_{\FU/\FV}
\end{alignat*}

\end{proof}

\begin{lem}
\item
\begin{enumerate}
\item
Let $X \in \sch$ and $\kappa \colon \FX:=X_{/X'} \to X$  the completion morphism of $X$ along a closed subscheme $X'\subset X$. Then  $\WCL_{\FX/X}\simeq 0$.
\item
If $\FX \xto{f} \FY$ is an \emph{adic} \'etale  morphism then $\WCL_{\FX/\FY}\simeq0$.
\end{enumerate}

\end{lem}

\begin{proof}
(1) By last proposition we can assume that $\kappa \colon \FX:=X_{/Z}=\spf(\widehat{A}) \to X=\spec(A)$ is the completion morphism of $X$ along a closed subscheme $Z= \spec(A/I)\subset X$. Put $\CI=I^{\tr}$. By Greenlees-May duality \cite[Theorem 0.1]{AJL3}, there is a $\D(\FX)$-isomorphism $\WCL_{\FX/X} \cong\LA_{\FX} (\CL_{\FX/X})$.

On the other hand, put $\FX= \dirlim {n} X_{n}$ with $X_n: = \spec(A/I^n), \forall n \in \NN$, such that the following diagrams are cartesian:
\[
\begin{tikzpicture}
\draw[white] (0cm,2.5cm) -- +(0: \linewidth)
node (21) [black, pos = 0.4] {$ \FX$}
node (22) [black, pos = 0.6] {$ X$};
\draw[white] (0cm,0.5cm) -- +(0: \linewidth)
node (11) [black, pos = 0.4] {$ X_n =\FX \times_X X_n$}
node (12) [black, pos = 0.6] {$X_n$.};
\draw [->] (21)  --   (22) node[above, midway, scale=0.75]{$\kappa$};
\draw [->] (11)  --   (12)node[above, midway, scale=0.75]{$1$};
\draw [right hook->] (11)  --   (21)node[left, midway, scale=0.75]{$i_n$};
\draw [right hook->] (12)  --   (22)node[right, midway, scale=0.75]{$j_n$};
\end{tikzpicture}
\]
Being $\kappa$  a flat morphism and  $X_n=\FX \times_X X_n$,  we have that the morphisms  $i_n^* \CL_{\FX/X}\xto{\sim
} \CL_{X_n/X_n}=0$ are isomorphisms \cite[Chapitre II, Corollaire 2.2.3]{I}.
Therefore, by \cite[\textbf{0}.13.2.3]{EGA31}  
\[
H^i(\LA_{\FX} (\CL_{\FX/X})) \simeq H^i(\invlim {n} i_n^* \CL_{\FX/X})\simeq \invlim {n} H^i( i_n^* \CL_{\FX/X})\simeq 0\]
for all $ i$, hence the conclusion.

Part (2) follows from analogous arguments.
\end{proof}

\begin{propo} \label{propcotcompl}
Given a commutative diagram in $\sfna$ 
\[
\begin{tikzpicture}[scale=0.75]
\draw[white] (0cm,2.5cm) -- +(0: \linewidth)
node (21) [black, pos = 0.4] {$X $}
node (22) [black, pos = 0.6] {$ Y$};
\draw[white] (0cm,0cm) -- +(0: \linewidth)
node (11) [black, pos = 0.4] {$ \FX$}
node (12) [black, pos = 0.6] {$\FY$};
\draw [->] (21)  --   (22) node[above, midway, scale=0.75]{$h$};
\draw [->] (11)  --   (12)node[above, midway, scale=0.75]{$\hh$};
\draw [->] (11)  --   (21)node[left, midway, scale=0.75]{$\kappa$};
\draw [->] (12)  --   (22)node[right, midway, scale=0.75]{$\kappa$};
\end{tikzpicture}
\]
with $\hh$ the completion morphism of a finite type morphism $h$ in $\scha$ along closed subschemes $X' \subset X, Y' \subset Y$ such that $h(X') \subset Y'$, and $\kappa$'s  completion morphisms, it holds that $\WCL_{\FX/\FY}\simeq
\WCL_{\FX/Y}\simeq\kappa^* \CL_{X/Y}$ and $\WCL_{\FX/\FY} \in \D^-_{\cc}(\FX)$.
\end{propo}

\begin{proof} Since $\CL_{X/Y} \in \D_{c}^{-}(X)$ \cite[Chapitre II, Corollaire 2.3.7]{I}, then $\kappa^*\CL_{X/Y} \in \D_{c}^{-}(\FX)$  and $\RLA_\FX k^* \CL_{X/Y} \cong  k^* \CL_{X/Y}$ by \ref{propihomologyloc}.(\ref{galacoh}). The result is immediate from the last lemma and the distinguished triangles associated to $h \circ \kappa$ and $\kappa \circ \hh$:
\begin{alignat*}{6}
\kappa^{*} \CL_{X/Y}& \to &\WCL_{\FX/Y}& \to& \WCL_{\FX/X}&\xto{+} \\
\RLA_\FX \LD \hh^{*} \WCL_{\FY/Y} &\to& \WCL_{\FX/Y} &\to& \WCL_{\FX/\FY} &\xto{+}. 
\end{alignat*}
\end{proof}
\begin{rem}
In the hypothesis of Proposition \ref{propcotcompl} it holds that $k^* \CL_{X/Y} \xto{\sim}  \WCL_{\FX/\FY}  $ is a $q$-flat resolution.
\end{rem}

\begin{propo} \label{coh}
Let $f \colon \FX \to \FY$ be a morphism of pseudo finite type in $\sfn$. 
Then $\WCL_{\FX/\FY} \in \D_{\cc}^{-}(\FX)$ and if $f \colon \FX= \spf(B) \to \FY = \spf(A)$ is in $\sfna$, then  $\h^{i}(\WCL_{\FX/\FY}) \simeq(\h^{i}(\WL_{B/A}))^{\tr}$.
\end{propo}

\begin{proof}
By Proposition \ref{localcotag} we can reduce to the affine case and consider that $f=\hh$, with $h$ a finite type morphism in $\scha$. Using Proposition \ref{propcotcompl} we can assume that $\FY=Y$ and therefore $f$ factorizes as  
\[
\begin{tikzpicture}
\draw[white] (0cm,1.5cm) -- +(0: \linewidth)
node (21) [black, pos = 0.2] {$\FX = \spf(B) $}
node (23) [black, pos = 0.5] {$ X=\spec(C)$}
node (24) [black, pos = 0.8] {$ Y=\spec(A)$};
\draw[white] (0cm,0cm) -- +(0: \linewidth)
node (12) [black, pos = 0.35] {$\spec(B)$};
\draw [->] (21)  --   (23) node[above, midway, scale=0.75]{$\kappa$};
\draw [->] (23)  --   (24)node[above, midway, scale=0.75]{$h$};
\draw [->] (21)  --   (12)node[above, midway, scale=0.75]{$\iota$};
\draw [->] (12)  --   (23);
\end{tikzpicture}
\]
with $B= \widehat{C}$ and $\iota$ is the canonical morphism.
Since $\kappa^*$ and $\iota^*$ are exact, it follows that

\begin{alignat*}{4}   
\h^{i}(\WCL_{\FX/\FY} ) &\simeq&\h^{i}( \kappa^* \CL_{X/Y} )& \qquad\text{Proposition \ref{propcotcompl}}\\
&\simeq &\h^{i}( \kappa^*   \widetilde{L_{C/A}})& \qquad\text{\cite[Chapitre II, Corollaire 2.3.7]{I}}\\
&\simeq&\quad  \iota^* ( \h^{i}(\widetilde{L_{C/A} \otimes_{C} B }))&\\
&\simeq&\quad \iota^* ((\h^{i}(L_{C/A} \otimes_{C} B ))^{\widetilde{\quad }}) &  \\
& \simeq& \quad \iota^* ((\h^{i}(\WL_{B/A}))^{\widetilde{\quad }})&\qquad \text{Corollary \ref{loc}}\\
& \simeq&  (\h^{i}(\WL_{B/A}))^{\tr}&\qquad \text{Prop. \ref{pftcoh}, \cite[(10.8.8)]{EGA1}}\\
\end{alignat*}
\end{proof}

\begin{cor} \label{h0}
Given $f \colon \FX \to \FY$  in $\sfn$, then $\h^{0}(\WCL_{\FX/\FY}) \simeq \om^{1}_{\FX/\FY}$.
\end{cor}
\begin{proof}
It follows from the last proposition and from Proposition \ref{h0afin}.
\end{proof}
\begin{parraf} \cite[3.4]{AJP1}
Given $\FX = \dirlim {n} X_n$ in $\sfn$, denote by $\A_{\mspace{1mu}\sphat\mspace{1mu}} (\FX)$  the full subcategory of  $\CO_{\FX}$-modules  $\CF$ such that 
\[\CF = \invlim {n \in \NN} (\CF \otimes_{\CO_{\FX}} \CO_{X_{n}}).\] 
For instance,  $\A_\cc(\FX)$ is a full subcategory of  $\A_{\mspace{1mu}\sphat\mspace{1mu}} (\FX)$ and if $\FX= \spf(A)$ and $M$ a $A$-module, $M^{\tr} \in \A_{\mspace{1mu}\sphat\mspace{1mu}} (\FX)$.
\end{parraf}

\begin{cor}
Let $f: \FX \to \FY $ be a morphism in $\sfn$. For every $\CF \in \A_{\mspace{1mu}\sphat\mspace{1mu}} (\FX)$ there is  a canonical isomorphism
\[
\ext^{0}_{\CO_{\FX}} (\WL_{\FX/\FY},\CF)  \overset{\sim}\lto  \Dercont_{\FY}(\CO_{\FX}, \CF)
\] 
\end{cor}
 \begin{proof}
 Use that $\Hom_{\CO_{\FX}}(\om^{1}_{\FX/\FY},\CF)\cong \Dercont_{\FY}(\CO_{\FX}, \CF) $
\cite[Theorem 3.5]{AJP1}.
\end{proof}

\begin{cor} \label{h0h1closed}
If $j \colon \FX \to \FY$ is a closed immersion in $\sfn$ given by an ideal $\CJ \subset \CO_{\FY}$, then $\WCL_{\FX/\FY} \simeq
\CL_{\FX/\FY}$. Therefore, $\h^{0}(\WCL_{\FX/\FY}) \simeq 0$ and $\h^{1}(\WCL_{\FX/\FY}) \simeq \CJ/\CJ^{2}$.
\end{cor}

\begin{proof} 
It suffices to apply Proposition \ref{coh} and Proposition \ref{cloim}.
\end{proof}

\begin{propo} \label{triancoh}
Let $f \colon \FX \to \FY$ and $g \colon \FY \to \FS$ be pseudo finite type morphisms in $\sfn$. Then there is a distinguished triangle in $\D_{\cc}^{-}(\FX)$
\[
\LD f^{*} \WCL_{\FY/\FS} \to \WCL_{\FX/\FS} \to \WCL_{\FX/\FY} \xto{+} 
\]
such that the maps in degree zero are the canonical ones.

\end{propo}
\begin{proof}
Combine Proposition \ref{triannocoh}, Proposition \ref{coh}  and  \ref{propihomologyloc}.(\ref{galacoh}). 
\end{proof}

\begin{parraf} \label{sef}

Let $f \colon \FX \to \FY$ and $g \colon\FY \to \FS$ be two pseudo finite type morphisms in $\sfn$.  By  Corollary \ref{h0},  the cohomology exact sequence  associated to the distinguished triangle (Proposition \ref{triancoh}) finishes, except to a functorial canonical isomorphism, in the \emph{complete} first fundamental exact sequence:
\[
\om^{1}_{\FY/\FS} \otimes \CO_{\FX} \to \om^{1}_{\FX/\FS} \to \om^{1}_{\FX/\FY} \to 0,
\] 
\cite[Proposition 3.10]{AJP1} (\cf\  \cite[Lemma 2.5.2]{LNS}).

If $f$ is a closed immersion given by an ideal $\CJ\subset\CO_{\FY}$, then the cohomology exact sequence  associated to the distinguished triangle finishes, except to a functorial canonical isomorphism, in the \emph{complete} second fundamental exact sequence (\cite[Proposition 3.13]{AJP1})
\[
\CJ/\CJ^2 \to \om^{1}_{\FY/\FS} \otimes \CO_{\FX} \to \om^{1}_{ \FX/\FS} \to  0.
\] 
\end{parraf}

\begin{propo}[Flat base-change]  \label{flatchange}
Consider a cartesian diagram  in $\sfn$ of pseudo finite type morphisms
\[
\begin{tikzpicture}[scale=0.75]
\draw[white] (0cm,2.5cm) -- +(0: \linewidth)
node (21) [black, pos = 0.4] {$\FX $}
node (22) [black, pos = 0.6] {$ \FY$};
\draw[white] (0cm,0.5cm) -- +(0: \linewidth)
node (11) [black, pos = 0.4] {$ \FX'$}
node (12) [black, pos = 0.6] {$\FY'$};
\draw [->] (21)  --   (22) node[above, midway, scale=0.75]{$f$};
\draw [->] (11)  --   (12)node[above, midway, scale=0.75]{$f'$};
\draw [->] (11)  --   (21)node[left, midway, scale=0.75]{$u$};
\draw [->] (12)  --   (22)node[right, midway, scale=0.75]{$v$};
\end{tikzpicture}
\]
such that 
$f$ or $v$ are flat.
Then the canonical maps are isomorphisms:
\begin{enumerate}
\item
$\LD u^{*}\WCL_{\FX/\FY} \to  \WCL_{\FX'/\FY'}$
\item
$\LD u^{*}\WCL_{\FX/\FY} \oplus \LD f'^{*} \WCL_{\FY'/\FY} \to  \WCL_{\FX'/\FY}$
\end{enumerate}

\end{propo}

\begin{proof}
We can reduce to the affine case, and assume that the diagram is a completion of a 
cartesian diagram   in $\scha$ with $f_1$ or $v_1$  flat morphisms:
\[
\begin{tikzpicture}
\draw[white] (0cm,3.5cm) -- +(0: \linewidth)
node (41) [black, pos = 0.3] {$\FX $}
node (42) [black, pos = 0.5] {$ \FY$};
\draw[white] (0cm,2.5cm) -- +(0: \linewidth)
node (31) [black, pos = 0.4] {$X $}
node (315) [black, pos = 0.5] {}
node (32) [black, pos = 0.6] {$ Y$};
\draw[white] (0cm,1.5cm) -- +(0: \linewidth)
node (21) [black, pos = 0.3] {$\FX' $}
node (215) [black, pos = 0.4] {}
node (22) [black, pos = 0.5] {$ \FY'$};
\draw[white] (0cm,0.5cm) -- +(0: \linewidth)
node (11) [black, pos = 0.4] {$ X'$}
node (12) [black, pos = 0.6] {$Y'$};
\draw [->] (41)  --   (42) node[above, midway, scale=0.75]{$f$};
\draw [-] (21)  --   (215);
\draw [->] (215)  --   (22)node[above, midway,scale=0.75]{$f'$};
\draw [->] (31)  --   (32) node[above, pos=0.25, scale=0.75]{$f_1$};
\draw [->] (11)  --   (12)node[above, midway, scale=0.75]{$f'_1$};
\draw [->] (21)  --   (41)node[left, pos=0.5, scale=0.75]{$u$};
\draw [-] (22)  --   (315)node[right, midway, scale=0.75]{};
\draw [->] (315)  --   (42)node[right, midway, scale=0.75]{$v$};
\draw [->] (11)  --   (31)node[right, pos=0.25, scale=0.75]{$u_1$};
\draw [->] (12)  --   (32)node[right, midway, scale=0.75]{$v_1$};
\draw [->] (41)  --   (31)node[above, pos=0.5, scale=0.75]{$\kappa$};
\draw [->] (42)  --   (32)node[above, pos=0.5, scale=0.75]{$\kappa$};
\draw [->] (21)  --   (11)node[above, pos=0.5, scale=0.75]{$\kappa$};
\draw [->] (22)  --   (12)node[above, pos=0.5, scale=0.75]{$\kappa$};
\end{tikzpicture}
\]
Since $\kappa^* \CL_{X/Y} \xto{\sim} \WCL_{\FX/\FY}$ is a $q$-flat resolution and $\kappa^*$ is exact,  from
the analogous result in $\sch$ \cite[Chapitre II, Proposition 2.2.3]{I},  it follows that:
\[
 \LD u^{*}\WCL_{\FX/\FY} = u^* \kappa^* \CL_{X/Y}=  \kappa^* u_1^*\CL_{X/Y} \simeq
\kappa^*\CL_{X'/Y'} \simeq \WCL_{\FX'/\FY'} 
\]

With similar considerations, the 
second part is consequence of \cite[Chapitre II, Proposition 2.2.3]{I}.

\end{proof}

\section{Deformation of formal schemes} \label{sec4} 

\begin{parraf}
\cite[Chapitre III, 1.1.10]{I}
Given a commutative diagram of ringed topoi
\begin{equation}
\begin{tikzpicture}
\draw[white] (0cm,2cm) -- +(0: \linewidth)
node (21) [black, pos = 0.425] {$ \FX$};
\draw[white] (0cm,0.5cm) -- +(0: \linewidth)
node (11) [black, pos = 0.425] {$ \FX'$}
node (12) [black, pos = 0.575] {$\FY$};
\draw [->] (11)  -- (12)node[above, midway, scale=0.75]{};
\draw [right hook->] (11)  --   (21)node[left, midway, scale=0.75]{$j$};
\draw [->] (21)  --   (12);
\end{tikzpicture}
\label{diagexte}
\end{equation}
where $j^{-1}(\CO_{\FX}) \to \CO_{\FX'}$ is a surjective map with square zero kernel we say that $\FX$  is a \emph{$\FY$-extension of $\FX'$}  (or $j \colon \FX' \inc \FX$ is a \emph{$\FY$-extension}). If $\CM$ is a $\CO_{\FX'}$-module such that $ \CM \simeq\Ker (j^{-1}(\CO_{\FX}) \to \CO_{\FX'})$, then  we say that $\FX$  is a \emph{$\FY$-extension of $\FX'$ by $\CM$}.
Notice that $\FX'_{\topo}= \FX_{\topo}$.

Next proposition ensures that the notion of  extension is stable in the category of locally noetherian formal schemes.
\end{parraf}

\begin{propo} \label{extenfs}
Let $(\FX', \CO_{\FX'})$ be in $\sfn$ and $\CM$ a coherent $\CO_{\FX'}$-module. Any extension as ringed topos of $\FX'$ by $\CM$  is in $\sfn$.
\end{propo}

\begin{proof}
It's well-known that any extension of $\FX'$ is a locally ringed space (\cf\   \cite[Chapitre III, 2.1.9.a)]{I}). Take $(\FX, \CO_{\FX})$ an extension of $\FX'$ by $\CM$ and  $\CI'$  an ideal of definition of $\FX'$.

Given $\FU' \subset \FX'$ an affine open formal subscheme, by  \cite[Corollary 3.1.8]{AJL2} $\h^{1}(\FU', \CM)=0$. We identify $\FU'$ with an open subset $\FU \subset \FX$. Then  we have a short exact sequence
\begin{equation} \label{sec}
0\to  M:=\ga(\FU, \CM) \to B:=\ga(\FU, \CO_{\FX}) \to B':=\ga(\FU, \CO_{\FX'})  \to 0 
\end{equation}
Since $B'$ and $M$ are noetherian we deduce that $B$ is also noetherian. 
Moreover, if $I'=\ga(\FU, \CI')$  and $I:=I'^{c}\subset B$, then $B$ is  $I$-adic. Indeed, since $M$ and $B'$ are a finitely generated $B$-modules and $I^{e}=I'$,  by $I$-adic completion of (\ref{sec}) we get a short exact sequence
\[
0 \to M    \to \widehat{B} \to B'  \to 0,
\]
from which we obtain that $B=\widehat{B}$.

Let $\CI\subset \CO_{\FX}$ be the ideal associated to the presheaf of ideals in $\FX$ locally defined by $ I=(\ga(\FU, \CI'))^{c}$, for all open subsets $\FU \subset \FX$ corresponding to open affine formal subschemes $ \spf(B') \subset \FX'$. 
From the construction it follows that $\CO_{\FX}= \invlim {n} \CO_{\FX}/\CI^{n+1}$.

We're going to shown that the inductive system $(\FX, \CO_{\FX}/\CI^{n+1})_{n \in \NN}$ verifies the hypothesis of \cite[10.6.3]{EGA1}, hence $(\FX, \CO_{\FX})= \dirlim {n } (\FX, \CO_{\FX}/\CI^{n+1})$ is a formal scheme.
For all $n \in \NN$ we have short exact sequences:
\[
0 \to (\CI^{n+1}+\CM)/\CI^{n+1} \to \CO_{\FX}/\CI^{n+1} \to \CO_{\FX'}/\CI'^{n+1} \to 0
\]
First, the ideal $ (\CI^{n+1}+\CM)/\CI^{n+1}$ is of zero square, so by \cite[Chapitre III, 2.1.9. b)]{I} $(\FX, \CO_{\FX}/\CI^{n+1})$ is a (usual) scheme. The canonical morphisms $\CO_{\FX}/\CI^{m+1}\to \CO_{\FX}/\CI^{n+1}$ are surjective, for all $m \ge n$. Last, the kernel  of $\CO_{\FX}/\CI^{n+2}\to \CO_{\FX}/\CI^{n+1}$ is nilpotent, $\forall n \in \NN$.
\end{proof}

\begin{propo} \label{propoext}
Let $f \colon \FX \to \FY$ and $g \colon \FY \to \FS$ be two pseudo finite type morphisms in $\sfn$ and $\CF \in \D_c(\FX)$. Then:
\begin{enumerate}
\item \label{propoext1}
$\ext^{i}_{\CO_{\FX}} (\WCL_{\FX/\FY}, \CF) \simeq \ext^{i}_{\CO_{\FX}} (\CL_{\FX/\FY}, \CF)$
\item \label{propoext2}
$\ext^{i}_{\CO_{\FX}} (\LD f^* \WCL_{\FY/\FS}, \CF) \simeq \ext^{i}_{\CO_{\FX}} (f^* \CL_{\FY/\FS}, \CF)$.
\end{enumerate}
\end{propo}

\begin{proof}
 By \ref{propihomologyloc}.(\ref{homlala}) and \ref{propihomologyloc}.(\ref{galacoh}) we have the canonical isomorphism
\begin{equation} \label{iso}
\R\shom ^{\bullet}(\RLA_{\FX}\CE, \CF)  \cong \R\shom ^{\bullet}(\RLA_{\FX}\CE, \RLA_{\FX}\CF)  \cong\R \shom ^{\bullet}(\CE, \CF).
\end{equation}
 On the other hand, there is a functorial isomorphism
\begin{alignat*}{3}
\R\shom ^{\bullet}(\LD f^* \WCL_{\FY/\FS}, \CF) &\simeq 
 \R\shom ^{\bullet}(\RLA_{\FX}\LD f^* \WCL_{\FY/\FS}, \CF)&\quad  \text{\ref{propihomologyloc}.(\ref{galacoh})}\\ 
 &\simeq \R\shom ^{\bullet}(\RLA_{\FX}f^* \CL_{\FY/\FS}, \CF)&\quad \text{Lemma \ref{lemderinver}}\\
&\simeq \R\shom ^{\bullet}(f^* \CL_{\FY/\FS}, \CF)& \text{(\ref{iso})}
\end{alignat*}
Since $ \R \Hom^{\bullet}(-, -) \simeq \R\Gamma\R \shom^{\bullet}(-,  -) $ \cite[Exercises 2.5.10 (b)]{L}, the result follows from the two above isomorphisms.
\end{proof}

\begin{rem}
Notice that by Proposition \ref{extenfs} and Proposition \ref{propoext} all the results of deformations of ringed topoi and morphisms of ringed topoi (\cf\ \cite[Chapitre III, \S 2]{I}) are true on $\sfn$. We only rewrite here the generic ones and left  for the reader more detailed formulations. 
\end{rem}

\begin{cor}[Deformation of formal schemes] \label{cordeformalsch}
Let  
\[
\begin{tikzpicture}
\draw[white] (0cm,2cm) -- +(0: \linewidth)
node (21) [black, pos = 0.35] {$ \FX$}
node (22) [black, pos = 0.5] {$ \FY$}
node (23) [black, pos = 0.65]{};
\draw[white] (0cm,0.5cm) -- +(0: \linewidth)
node (11) [black, pos = 0.35] {$ \FX'$}
node (12) [black, pos = 0.5] {$\FY'$}
node (13) [black, pos = 0.65]{$\FS$};
\draw [dashed,->] (21)  --   (22) node[above, midway, scale=0.75]{$f$};
\draw [->] (11)  -- (12)node[above, midway, scale=0.75]{$f'$};
\draw [->] (12)  --   (13)node[above, midway, scale=0.75]{};
\draw [dashed,right hook->] (11)  --   (21)node[left, midway, scale=0.75]{$i$};
\draw [right hook->] (12)  --   (22)node[right, midway, scale=0.75]{$j$};
\draw [->] (22)  --   (13);
\end{tikzpicture}
\]
be  a commutative diagram of pseudo finite-type morphisms in $\sfn$  where 
 $j \colon \FY'  \inc \FY$ is a $\FS$-extension by a coherent $\CO_{\FY'}$-module $\CN $.

\begin{enumerate}
\item
(\emph{Existence}) Given a morphism of coherent $\CO_{\FX'}$-modules $v \colon f'^{*} \CN \to \CM$, there is an element $c(f',j,v) \in \ext^{2}_{\CO_{\FX'}}(\WCL_{\FX'/\FY'}, \CM)$  such that $c(f',j,v)=0$ if, and only if, there exists a lifting morphism $f \colon\FX \to \FY$ in $\sfn$ such that $i \colon \FX'  \inc \FX$ is  a $\FS$-extension by  $\CM$ and the above diagram is commutative.

\item 
(\emph{Uniqueness})
In the case that there exists a lifting morphism, 
the set of  isomorphism classes of such lifting morphisms  is an affine  space over $\ext^{1}_{\CO_{\FX'}}(\WCL_{\FX'/\FY'}, \CM)$.
\end{enumerate}
\end{cor}

\begin{proof}
Combine  Proposition \ref{propoext}.(\ref{propoext1}) and the analogous result of ringed topoi \cite[Chapitre III, Th\`eoreme 2.1.7]{I}.
\end{proof}

\begin{cor} [Deformation of morphisms of formal schemes] \label{cordemorfisformalsch}
Let 
\[
\begin{tikzpicture}
\draw[white] (0cm,2cm) -- +(0: \linewidth)
node (21) [black, pos = 0.25] {$ \FX$}
node (22) [black, pos = 0.4] {$ \FY$}
node (23) [black, pos = 0.55]{$\FW$};
\draw[white] (0cm,0.5cm) -- +(0: \linewidth)
node (11) [black, pos = 0.25] {$ \FX'$}
node (12) [black, pos = 0.4] {$\FY'$}
node (13) [black, pos = 0.55]{$\FW'$}
node (14) [black, pos = 0.7]{$\FS$};
\draw [dashed,->] (21)  --   (22) node[above, midway, scale=0.75]{$f$};
\draw [->] (22)  --   (23)node[above, midway, scale=0.75]{};
\draw [->] (11)  --   (12)node[above, midway, scale=0.75]{$f'$};
\draw [->] (12)  --   (13)node[above, midway, scale=0.75]{$g'$};
\draw [->] (13)  --   (14)node[above, midway, scale=0.75]{$q'$};
\draw [right hook->] (11)  --   (21)node[left, midway, scale=0.75]{$i$};
\draw [right hook->] (12)  --   (22)node[right, midway, scale=0.75]{$j$};
\draw [right hook->] (13)  --   (23)node[right, midway, scale=0.75]{$k$};
\draw [->] (23)  --   (14);
\end{tikzpicture}
\]
be  a commutative diagram of pseudo finite-type morphisms in $\sfn$  where 
$i \colon \FX'  \inc \FX$ is a $\FW$-extension by a coherent $\CO_{\FX'}$-module  $\CM$,  $j \colon \FY'  \inc \FY$ is a $\FW$-extension by a coherent $\CO_{\FY'}$-module   $\CN$ and $k \colon \FW'  \inc \FW$ is  a $\FS$-extension by a coherent $\CO_{\FW'}$-module  $\CL$. Take $v \colon g'^* \CL \to \CN$ and $w \colon (g' \circ f')^* \CL \to \CM$ two morphisms.  
\begin{enumerate}
\item
(\emph{Existence}) 
Given $u \colon f'^{*} \CN \to \CM$ a morphism of  $\CO_{\FX'}$-modules with $w = u \circ f'^* v$, there is  $c(f',i,j,u) \in \ext^{1}_{\CO_{\FX'}}(\LD f'^*\WCL_{\FY'/\FW'}, \CM)$ whose vanishing is equivalent to the existence of a lifting morphism $f \colon\FX \to \FY$ in $\sfn$ such that the above diagram is commutative.
\item
(\emph{Uniqueness}) 
In the case of existence,  the set of   lifting morphisms $\FX \to \FY$  is an affine  space over $\ext^{0}_{\CO_{\FX'}}(\LD f'^*\WCL_{\FY'/\FW'}, \CM)$.
\end{enumerate}
\end{cor}

\begin{proof}
Use  Proposition \ref{propoext}.(\ref{propoext2})  and the analogous result of ringed topoi \cite[Chapitre III, Proposition 2.2.4]{I}.
\end{proof}

\begin{rem}
We point out that the deformation theory presented here generalizes the obstruction results obtained in \cite{P} through the complete differential module and under smoothness  and separation hypothesis. 

\end{rem}
\section{Cotangent complex and infinitesimal conditions}\label{sec5}

\begin{propo} \label{cotanliso}
Let  $f \colon \FX \to \FY$ be a pseudo finite type morphism in $\sfn$. Then $f$ is smooth if, and only if, $\WCL_{\FX/\FY} \simeq
\om^{1}_{\FX/\FY}[0]$ and $\om^{1}_{\FX/\FY}$ is a finite-rank locally free $\CO_{\FX}$-module. 
\end{propo}

\begin{proof}
By \cite[Proposition 4.1]{AJP1} and  Proposition \ref{localcotag} we may assume that $f$ is in $\sfna$.  
The ``only if" is  consequence of \cite[Proposition 4.8]{AJP1},  Proposition \ref{coh} and  Proposition \ref{smo}.

For the converse we have to prove that if $Z$ is an affine $\FY$-scheme, $T\inc Z$ is a closed subscheme given by a square zero ideal $\CJ \subset \CO_{Z}$, then any  $\FY$-morphism $u\colon T \to \FX$ admits a lifting $Z \to \FX$. From the hypothesis we deduce that
\[
\ext^{1}_{\CO_{T}}(\LD u^{*} \WCL_{\FX/\FY}, \CJ)=\ext^{1}_{\CO_{T}}( u^{*}  \om^{1}_{\FX/\FY}[0], \CJ)=0
\]
since $\om^{1}_{\FX/\FY}[0]$ is a finite-rank locally free $\CO_{\FX}$-module an $T$ is an affine scheme.
The  conclusion follows from Corollary \ref{cordemorfisformalsch}.
\end{proof}

\begin{rem}
A well-known fact in the case of usual schemes, is that the above proposition leads to simpler versions of deformation results.  That is, whenever the involved morphism is smooth,  in the  Ext-groups the cotangent complex is replaced by the differential module. 
However, frequently  arise situations when the base morphism  is not smooth. Proposition \ref{cotanliso} provides  useful consequences of Corollaries \ref{cordeformalsch} and  \ref{cordemorfisformalsch} in the case of base scheme morphisms $X' \to Y$ where $X'$ is not smooth but closed embeddable in a smooth scheme $X$.  We get a smooth formal scheme $\FX:= X_{/X'}= \dirlim {n \in \NN}X_n$ such that $X_0=X'$ . Then we can obtain deformation conditions in terms of $ \om^{1}_{\FX/Y}[0]$ instead of $\WCL_{\FX/Y}$,  which is a more untreatable object. We give in \ref{apt} an example of this application.  
\end{rem}

\begin{parraf}\label{apt}
Let  $X' \xto{f'}Y'$ be a $Z$-morphism of schemes where $X'$ is  affine,  $j \colon Y' \inc Y $  a closed  immersion in a smooth $Z$-scheme  $Y$ and $\FY: =Y_{/Y'}$. The induced morphism $\FY \to Z$ is a smooth morphism \cite[Proposition 4.5]{AJP1} and $ \WCL_{\FY/Z} \simeq \om^{1}_{\FY/Z}[0]$ (Proposition \ref{cotanliso}). Given $i \colon X'  \inc X$  a $Z$-extension by a coherent $\CO_{X'}$-module  $\CM$, it holds that
\[
\ext^{1}_{\CO_{X'}}(\LD (j \circ f')^{*} \WCL_{\FY/Z}, \CM)=\ext^{1}_{\CO_{X'}}( (j \circ f')^{*}  \om^{1}_{\FY/Z}[0], \CM)=0
\]
Therefore by Corollary \ref{cordemorfisformalsch} there exists a morphism $X \xto{f} \FY$ such that the following diagram is commutative. \[
\begin{tikzpicture}
\draw[white] (0cm,2cm) -- +(0: \linewidth)
node (21) [black, pos = 0.25] {$ X$}
node (22) [black, pos = 0.4] {$ \FY=Y_{/Y'}$};
\draw[white] (0cm,0.5cm) -- +(0: \linewidth)
node (11) [black, pos = 0.25] {$ X'$}
node (12) [black, pos = 0.4] {$Y'$}
node (13) [black, pos = 0.55]{$Z$};
\draw [dashed,->] (21)  --   (22) node[above, midway, scale=0.75]{$f$};
\draw [->] (11)  --   (12)node[above, midway, scale=0.75]{$f'$};
\draw [->] (12)  --   (13)node[above, midway, scale=0.75]{$g'$};
\draw [right hook->] (11)  --   (21)node[left, midway, scale=0.75]{$i$};
\draw [right hook->] (12)  --   (22)node[right, midway, scale=0.75]{$j$};
\draw [->] (22)  --   (13);
\end{tikzpicture}
\]

\end{parraf}

\begin{cor} \label{cotanetale}
Let  $f \colon \FX \to \FY$ be a pseudo finite type morphism in $\sfn$. Then $f$ is \'etale if, and only if, $\WCL_{\FX/\FY} \simeq
0$.
\end{cor}

\begin{proof}

It is immediate from the last proposition and \cite[Proposition 4.6]{AJP1}.

\end{proof}

\begin{cor} \label{corclosedsmooth}
Let $j \colon \FX' \inc \FX$ be a closed immersion given by an  ideal $\CJ\subset \CO_{\FX}$ and $f \colon \FX \to \FY$ a smooth morphism. Then
\[
\tau^{\scriptscriptstyle{\ge -1}}(\WCL_{\FX'/\FY})  \simeq
(0 \to \CJ/\CJ^{2} \xto{\hd} j^*\om^{1}_{\FX/\FY}\to 0)
\]
with $\hd$ given by $\hd_{\FX/\FY}$.
\end{cor}

\begin{proof}
By Proposition \ref{cotanliso} the distinguished triangle  associated to $ f \circ j$ (Proposition \ref{triancoh}) is given by
\[
j^{*} \om^1_{\FX/\FY} \to \WCL_{\FX'/\FY} \to \WCL_{\FX'/\FX} \xto{+}. 
\]
Analogously to \ref{cotsef2}, using Corollary \ref{h0h1closed}, it is easy to show that the following diagram is commutative (\cf\ \cite[Chapitre III, Proposition 1.2.9]{I}):
\[
\begin{tikzpicture}
\draw[white] (0cm,2cm) -- +(0: \linewidth)
node (21) [black, pos = 0.35] {$ \h^{1}(\WCL_{\FX'/\FX})$}
node (22) [black, pos = 0.65] {$ \h^{0}(\LD j^{*} \WCL_{\FX/\FY})$};
\draw[white] (0cm,0.5cm) -- +(0: \linewidth)
node (11) [black, pos = 0.35] {$\CJ/\CJ^2$}
node (12) [black, pos = 0.65] {$j^*\om^{1}_{\FX/\FY}$};
\draw [->] (21)  --   (22);
\draw [->] (11)  --   (12)node[above, midway, scale=0.75]{$-\hd$};
\draw [->] (21)  --   (11)node[left, midway, scale=0.75]{$\simeq$} ;
\draw [->] (22)  --   (12)node[right, midway, scale=0.75]{$\simeq$};
\end{tikzpicture}
\]
with $\hd$ induced by $\hd_{\FX/\FY}$.
Then by truncation  we obtain a distinguished triangle 
\[
j^{*} \om^1_{\FX/\FY} \to \tau^{\scriptscriptstyle{\ge -1}}(\WCL_{\FX'/\FY}) \to \CJ/\CJ^{2}[1] \xto{+} 
\]
such that the connecting morphism is $-\hd$, therefore the conclusion.
\end{proof}

\begin{defn}
Let $j \colon  \FX' \inc \FX$ be a closed immersion in $\sfn$ given by an ideal $\CJ \subset \CO_{\FX}$. We say that $j$ is a \emph{regular closed immersion} if $\CJ$ is regular \cite[16.9.1]{EGA44}.

\end{defn}

\begin{parraf} \label{proporegularclosed}
If  $j \colon \FX' \inc \FX$ is a regular closed immersion with ideal $\CJ\subset \CO_{\FX}$, then $\CJ/\CJ^2 $  is a locally free $\CO_{\FX}/\CJ$-module and $\WCL_{\FX'/\FX} \simeq\CJ/\CJ^2 [1]$.
Indeed, the former follows from \cite[16.9.4]{EGA44} and the second it is immediate from \cite[Chapitre III, Proposition 3.2.4]{I},
since $\WCL_{\FX'/\FX} \simeq
\CL_{\FX'/\FX}$ (Corollary \ref{h0h1closed}).
\end{parraf}

\begin{propo} \label{complinters}
Let $j \colon \FX' \inc \FX$ be a regular closed immersion given by an ideal $\CJ\subset \CO_{\FX}$ and $f \colon \FX \to \FY$ a smooth morphism. Then $\WCL_{\FX'/\FY}$ is isomorphic to a complex of locally free $\CO_{\FX'}$-modules concentrated in degree $[-1,0]$:
\[
\WCL_{\FX'/\FY}  \simeq
(0 \to \CJ/\CJ^{2} \to j^*\om^{1}_{\FX/\FY}\to 0).
\]
\end{propo}

\begin{proof}
Apply Corollary \ref{corclosedsmooth} and  \ref{proporegularclosed} to the distinguished triangle associated to 
$ \FX' \inc \FX \to \FY$.
\end{proof}

\begin{rem}
Most of the proofs in this section are a rewritten, \emph{mutatis mutandis}, of the analogous ones in \cite[Chapitre III, \S3]{I}, but because of its dependence on the results of this paper, we included it here. 
\end{rem}

\end{document}